\documentclass[12pt]{article}

\usepackage{latexsym,amsmath,amscd,amssymb,framed,graphics,mathrsfs}
\usepackage{enumerate}

\usepackage{graphicx}

\usepackage[colorlinks]{hyperref}
\usepackage{url}
\usepackage{colordvi}
\usepackage{color}
\usepackage{bbm}
\usepackage{cancel}

\usepackage{soul}
\setstcolor{blue}

\usepackage[all]{xy}

\makeatletter

\@addtoreset{figure}{section}
\def\thefigure{\thesection.\@arabic\c@figure}
\def\fps@figure{h, t}
\@addtoreset{table}{bsection}
\def\thetable{\thesection.\@arabic\c@table}
\def\fps@table{h, t}
\@addtoreset{equation}{section}
\def\theequation{\thesection.\arabic{equation}}
\makeatother

\newenvironment{proof}{\paragraph{Proof:}}{\hfill$\square$}
\newtheorem{theorem}{Theorem}[section]

\newtheorem{lemma}[theorem]{Lemma}
\newtheorem{remark}[theorem]{Remark}
\newtheorem{proposition}[theorem]{Proposition}

\begin{document}

\title{Variational integrators for anelastic and pseudo-incompressible flows}
\author{Werner Bauer$^{1,2}$ and Fran\c{c}ois Gay-Balmaz$^2$}
\addtocounter{footnote}{1} 
\footnotetext{Imperial College London, Department of Mathematics. 
\texttt{w.bauer@imperial.ac.uk}}
\addtocounter{footnote}{1} 
\footnotetext{CNRS and \'Ecole Normale Sup\'erieure, Laboratoire de M\'et\'eorologie Dynamique, Paris, France. 
\texttt{gaybalma@lmd.ens.fr}}

\date{\fontsize{11pt}{12pt}\selectfont To Darryl Holm, on the occasion of his 70${}^{th}$ birthday}
\maketitle

\begin{abstract}
  The anelastic and pseudo-incompressible equations are two well-known soundproof approximations 
  of compressible flows useful for both theoretical and numerical analysis in meteorology, atmospheric 
  science, and ocean studies. 
  In this paper, we derive and test structure-preserving numerical schemes for these two systems.
  The derivations are based on a discrete version of the Euler-Poincar\'e variational method. 
  This approach relies on a finite dimensional approximation of the (Lie) group of diffeomorphisms that 
  preserve {\em weighted-}volume forms. These weights describe the background stratification of the 
  fluid and correspond to the weighed velocity fields for anelastic and pseudo-incompressible approximations.
  In particular, we identify to these discrete Lie group configurations the 
  associated Lie algebras such that elements of the latter correspond to weighted velocity fields that 
  satisfy the divergence-free conditions for both systems. Defining discrete Lagrangians in terms of these 
  Lie algebras, the discrete equations follow by means of variational principles.  
  Descending from variational principles, the schemes exhibit further a 
  discrete version of Kelvin circulation theorem, are applicable to irregular meshes, and 
  show excellent long term energy behavior. We illustrate the properties of the schemes by performing preliminary test cases. 
\end{abstract}

\section{Introduction}
 
 Numerical simulations of atmosphere and ocean on the global scale are of high 
 importance in the field of Geophysical Fluid Dynamics (GFD).
 The dynamics of these systems are frequently modeled by the full Euler equations using explicit time 
 integration schemes (see, e.g., \cite{DurranGFD1999}). These simulations are however computationally very expensive. 
 Besides highly resolved meshes to capture important small scale features, the fast 
 traveling sound waves have to be resolved too, by very small time step sizes,
 in order to guarantee stable simulations \cite{DurranGFD1999}. As these sound waves are assumed to be 
 negligible in atmospheric flows, soundproof models, in which these fast waves are 
 filtered out, are a viable option that permits to increase the time step sizes and 
 hence to speed up calculations significantly.

 Frequently applied soundproof models are the Boussinesq, anelastic, and pseudo-incompressible 
 approximations of the full Euler equations \cite{Du1989,LiHe1982,OgPh1962}. 
 There exist elaborated discretizations of these equations
 in literature. However, these discretizations often do not take into account 
 the underlying geometrical structure of the equations. This may result in a lack of 
 conserving mass, momentum, energy, or to the fact that the Helmholtz decomposition 
 of vector fields or the Kelvin-Noether circulation theorem are not satisfied. 
 Structure-preserving schemes descending from Euler-Poincar\'e variational methods \cite{PaMuToKaMaDe2010}, \cite{GaMuPaMaDe2010}, \cite{DeGaGBZe2014} conserve these 
 quantities, as they arise from a Lagrangian formulation, in which these conserved quantities
 are given by invariants of the Lagrangian under symmetries, \cite{HoMaRa1998}.

 With this paper, we contribute to develop variational integrators in the area of GFD 
 by including the anelastic and pseudo-incompressible schemes into the variational discretization
 framework developed by \cite{PaMuToKaMaDe2010}. To use this framework, we first have to describe 
 these approximations of the Euler equations in terms of the Euler-Poincar\'e variational 
 method \cite{HoMaRa1998}. The central idea is to use volume forms that are weighted by the corresponding background 
 stratifications such that they match the divergence-free conditions of the correspondingly weighed velocity
 fields associated to either the anelastic or the pseudo-incompressible approximations. 
 This will allow us to identify for these approximations the appropriate Lie group configuration 
 with corresponding Lie algebras. Using the latter to define appropriate Lagrangians,
 the equations of motion follow by Hamilton's variational principle of stationary action.
 
 The definition of appropriate discrete diffeomorphism groups will be based on the
 idea to use weighted meshes that provide discrete counterparts of the weighted volume forms.
 The corresponding discrete Lie algebras will incorporate the required divergence-free conditions 
 on the weighed velocity fields. Defining appropriate weighted pairings required to derive the 
 functional derivatives of the discrete Lagrangians, the flat operator introduced in \cite{PaMuToKaMaDe2010}
 is directly applicable and we can thus avoid to discuss this otherwise delicate issue.
 Mimicking the continuous theory, the discretizations of anelastic and pseudo-incompressible 
 equations follow by variations of appropriate discrete Lagrangians.

 We structure the paper as follows. 
 In Section \ref{sec_anelastic_pseudoincomp_systems} we recall the standard formulations
 of Boussinesq, anelastic, and pseudo-incompressible approximations of the Euler equations 
 for perfect fluids. 
 In Section \ref{VF} we show that these equations follow from the Euler-Poincar\'e variational principle, for appropriate 
 Langrangians. 
 In Section \ref{VD} we recall the variational discretization
 framework introduced by \cite{PaMuToKaMaDe2010} and extend it to suit 
 anelastic and pseudo-incompressible equations. The corresponding 
 discretizations on 2D simplicial meshes are presented in Section \ref{sec_numerical_schemes}, and preliminary numerical tests are performed in \ref{sec_numerical_tests}.
 In Section \ref{sec_conclusions} we draw conclusions and provide an outlook.

\section{Anelastic and pseudo-incompressible systems}
\label{sec_anelastic_pseudoincomp_systems}

In this section we review the three approximations of the Euler equations of a perfect gas that will be the subject of this paper, namely, 
the Boussinesq, the anelastic, and the pseudo-incompressible approximations
(see, e.g., \cite{DurranGFD1999} for more details).

The Euler equations for the inviscid isentropic motion of a perfect gas can be expressed in the form
\begin{equation}\label{atmosph_motion} 
\partial _t \mathbf{u} + \mathbf{u} \cdot \nabla \mathbf{u} +\frac{1}{\rho }\nabla p =- g \mathbf{z} , \quad \partial _t \rho + \operatorname{div}( \rho \mathbf{u} )=0, \quad \partial _t \theta + \mathbf{u} \cdot \nabla \theta =0,
\end{equation} 
where $ \mathbf{u} $ is the three-dimensional velocity vector, $ \rho$ is the mass density, $p$ is the pressure, $g$ is the gravitational acceleration, 
$ \mathbf{z} $ is the unit vector directed opposite to the gravitational force. 
The variable $ \theta$ is the potential temperature, defined by $ \theta =T/ \pi $, 
in which $T$ is the temperature and $ \pi $ is the Exner pressure
\[
\pi =(p/p_0)^{R/c_p},
\]
with $R$ the gas constant for dry air, $c_p$ the specific heat at constant pressure, 
and $p_0$ a constant reference pressure. Using the equation of state for a perfect gas, $p= \rho RT$, we have the relation
\[
\frac{1}{\rho } \nabla p= c_p\theta \nabla \pi .
\]
The equations \eqref{atmosph_motion} correspond to conservation of momentum, mass, and entropy, respectively.

Let us write
\[
\theta  (x,y,z,t)=\bar\theta  (z)+ \theta  '(x,y,z,t), \quad \pi (x,y,z,t)=\bar \pi  (z)+ \pi '(x,y,z,t),
\]
in which $\bar\theta  (z)$ and $\bar \pi (z)$ characterize a vertically varying reference state in hydrostatic balance, that is,
\begin{equation}\label{HB} 
c _p \bar \theta \frac{d \bar \pi }{d z}=- g.
\end{equation} 
In terms of the perturbations $ \theta '$ and $ \pi '$, the equations \eqref{atmosph_motion} can be equivalently written as
\begin{equation}\label{equiv} 
\partial _t \mathbf{u} + \mathbf{u} \cdot \nabla \mathbf{u}+ c _p \theta \nabla \pi '= g \frac{\theta '}{\bar \theta } \mathbf{z} , \quad \partial _t \rho + \operatorname{div}( \rho \mathbf{u} )=0, \quad \partial _t \theta + \mathbf{u} \cdot \nabla \theta =0.
\end{equation} 

We introduce now three frequently applied approximations to these equations.

\paragraph{Boussinesq approximation.} This approximation is obtained by assuming a nondivergent flow 
and by neglecting the variations in potential temperature except in the leading-order contribution to 
the buoyancy. We thus get, from \eqref{equiv}, the system
\[
\partial _t \mathbf{u} + \mathbf{u} \cdot \nabla \mathbf{u}+ c _p \theta_0 \nabla \pi '= g \frac{\theta '}{ \theta _0} \mathbf{z} , 
\quad \operatorname{div}\mathbf{u} =0, \quad \partial _t \theta + \mathbf{u} \cdot \nabla \theta =0,
\]
in which $\theta_0$ is a constant reference potential temperature. These equations can equivalently be written as
\[
\partial _t \mathbf{u} + \mathbf{u} \cdot \nabla \mathbf{u}+ \nabla P'_{\rm b}= b'\mathbf{z} , \quad \operatorname{div}\mathbf{u} =0, \quad \partial _t b' + \mathbf{u} \cdot \nabla b'+N ^2 w =0,
\]
where $P'_{\rm b}= c_p \theta _0 \pi '$, $N^2 = \frac{g}{ \theta _0} \partial _z \bar \theta $  is the Brunt-V\"ais\"al\"a frequency, and $ b '=g\frac{ \theta '}{\theta _0} $ is the buoyancy. 
Making use of the full buoyancy $b= g \frac{\theta }{\theta _0}=  g \frac{\bar \theta+\theta ' }{\theta _0}$, we can write the system as
\begin{equation}\label{Boussinesq_standard_form} 
\partial _t \mathbf{u} + \mathbf{u} \cdot \nabla \mathbf{u}+ \nabla P_{\rm b} = b\mathbf{z} , \quad \operatorname{div}\mathbf{u} =0, \quad \partial _t b + \mathbf{u} \cdot \nabla b =0,
\end{equation} 
where $P_{\rm b} := P'_{\rm b} + \frac{g}{\theta_0} \int_0^z \bar \theta(z) dz$.

The total energy is conserved since the energy density $E=\frac{1}{2} | \mathbf{u} |^2 -bz=\frac{1}{2} |\mathbf{u} | ^2 - g \frac{\theta }{\theta _0}z$ verifies the continuity equation
\begin{equation}\label{energy_cons} 
\partial _t E+ \operatorname{div}((E+ P_{\rm b}) \mathbf{u} ) =0 .
\end{equation} 
The requirement for nondivergent flow is easily justified only for liquids, and the errors incurred approximating the true mass 
conservation relation by $\operatorname{div}\mathbf{u} =0$ can be quite large in stratified compressible flows. In this case, 
the anelastic and pseudo-incompressible models have to be considered, which better approximate the true mass continuity equation.

\paragraph{Anelastic approximation.} The anelastic system approximates the continuity equation as 
\[
\operatorname{div}(\bar \rho\mathbf{u} )=0, 
\]
where $\bar \rho (z)$ is the vertically varying density of the reference state.

In the original anelastic system presented by \cite{OgPh1962}, the reference state is isentropic so that $\bar \theta (z) = \theta _0$ is constant, which results in the approximation
\begin{equation}\label{anelastic0} 
\partial _t \mathbf{u} + \mathbf{u} \cdot \nabla \mathbf{u}+ c _p \theta _0\nabla \pi '= g \frac{\theta '}{\theta _0} \mathbf{z} , \quad \operatorname{div}(\bar \rho\mathbf{u} )=0, \quad \partial _t \theta + \mathbf{u} \cdot \nabla \theta =0.
\end{equation} 
The energy density can be written as $E= \bar \rho \left( \frac{1}{2} |\mathbf{u} | ^2 - g \frac{\theta }{\theta _0}z\right) = \bar \rho\left(  \frac{1}{2} | \mathbf{u} | ^2 + c _p \bar \pi \theta \right)$, where $\bar \pi (z)= - \frac{g}{c_p \theta _0}z$ verifies the hydrostatic balance \eqref{HB} for $\bar \theta (z)= \theta _0$. The total energy is conserved since $E$ verifies the continuity equation
\[
\partial _t E+ \operatorname{div}( (E+ P_{\rm a_0}) \mathbf{u} )=0
\]
with $P_{\rm a_0} := \bar{ \rho }(c_p \theta _0 \pi '+ gz)$
  
In the subsequent work \cite{WiOg1972}, the reference potential temperature $\bar \theta $ was allowed to vary in the vertical, leading to the momentum equation
\[
\partial _t \mathbf{u} + \mathbf{u} \cdot \nabla \mathbf{u}+ c _p\bar \theta \nabla \pi '= g \frac{\theta '}{\bar \theta } \mathbf{z}.
\]
The resulting system is however not energy conservative. In order to restore energy conservation, \cite{LiHe1982} considered the approximate momentum equation
\begin{equation}\label{anelastic} 
\partial _t \mathbf{u} + \mathbf{u} \cdot \nabla \mathbf{u}+ \nabla(c _p \bar \theta \pi ')= g \frac{\theta '}{\bar \theta } \mathbf{z} , \quad \operatorname{div}(\bar \rho\mathbf{u} )=0, \quad \partial _t \theta + \mathbf{u} \cdot \nabla \theta =0.
\end{equation} 
In this case, the energy density $E= \bar \rho\left(  \frac{1}{2} | \mathbf{u} | ^2 + c _p \bar \pi \theta \right)$, where $\bar\pi (z)$ is such that $ c_p \frac{\partial \bar \pi }{\partial z}=- \frac{g}{ \bar\theta }$, satisfies the continuity equation
\[
\partial _t E+ \operatorname{div}( (E+P_{\rm a})\mathbf{u} )=0, 
\]
with $P_{\rm a} := \bar\rho (c_p \bar \theta  \pi '+ gz)$.

\paragraph{Pseudo-incompressible approximation.} To obtain this approximation developed in \cite{Du1989}, one defines the pseudo-density $ \rho ^\ast = \bar\rho \bar\theta /\theta $ and enforces mass conservation with respect to $ \rho ^\ast $ as $ \partial_t  \rho ^\ast+ \operatorname{div}(\rho ^\ast \mathbf{u} )=0$. When combined with $ \partial _t\theta + \mathbf{u} \cdot \nabla \theta =0$, it yields $ \operatorname{div}( \bar \rho \bar \theta \mathbf{u} )=0$. These last two equations can be used with the momentum equation in \eqref{equiv} to yield the energy conservative system
\begin{equation}\label{pseudo_incompressible} 
\partial _t \mathbf{u} + \mathbf{u} \cdot \nabla \mathbf{u} + c _p \theta \nabla \pi ' = g \frac{\theta'}{\bar\theta} \mathbf{z} , \quad \operatorname{div}( \bar \rho \bar \theta \mathbf{u} )=0, \quad \partial _t \theta  + \mathbf{u} \cdot \nabla \theta =0.
\end{equation} 
We note that the balance of momentum is equivalently written as $\partial _t \mathbf{u} + \mathbf{u} \cdot \nabla \mathbf{u} + c _p \theta \nabla \pi = - g \mathbf{z}$, where $ \pi =\bar \pi +\pi '$, with $ c_p \frac{\partial \bar \pi }{\partial z}=- \frac{g}{ \bar\theta }$.
The energy density $E=\rho ^\ast \left(  \frac{1}{2} | \mathbf{u} | ^2 + gz \right)$ verifies the continuity equation
\[
\partial _t E+ \operatorname{div}((E+P_{\rm pi}) \mathbf{u} )=0
\]
for $P_{\rm pi}: = c_ p  \rho^\ast \theta\pi$.

\section{Variational formulation}\label{VF}

We shall now formulate the anelastic and pseudo-incompressible equations in Euler-Poincar\'e variational form. Euler-Poincar\'e variational principles are Eulerian versions of the classical Hamilton principle of critical action. We refer to \cite{HoMaRa1998} for the general Euler-Poincar\'e theory based on Lagrangian reduction and for several applications in fluid dynamics. 
An Euler-Poincar\'e formulation for anelastic systems was given in \cite{CoHo2013}. We shall develop below a slightly different 
Euler-Poincar\'e approach, well-suited for the variational discretization, by putting the emphasis on the underlying Lie group of diffeomorphisms associated to these systems.

%\todo{\color{magenta}FGB: From here on, I revised the text and added details and proof of all statements. I also assumed that $ \mathcal{D} $ has a boundary because this is what we do in the numerical schemes. This is much more rigorous now.}

As we have recalled above, the anelastic and pseudo-incompressible equations are based on a constraint of the following type on the fluid velocity $\mathbf{u} (t,\mathbf{x} )$:
\[
\operatorname{div}( \bar \sigma \mathbf{u} ) =0,
\] 
for a given strictly positive function $\bar \sigma ( \mathbf{x} )>0$ on the fluid domain $ \mathcal{D} $.

We assume that the fluid domain $ \mathcal{D} $ is a compact, connected, orientable manifold with smooth boundary 
$ \partial \mathcal{D} $. In our examples, $\mathcal{D} $  is a 2D domain in the vertical plane $ \mathbb{R}  ^2 \ni \mathbf{x} = (x,z)$ 
or a 3D domain in $ \mathbb{R}  ^3 \ni \mathbf{x} =(x,y,z)$.

We fix a volume form $ \mu $ on $ \mathcal{D} $, i.e., an $n$-form, $n= \operatorname{dim} \mathcal{D} $, with $ \mu (\mathbf{x} )\neq 0$, for all $ \mathbf{x} \in\mathcal{D} $ . If $ \mathcal{D}$ is a domain in $ \mathbb{R}  ^3\ni (x,y,z )$, one can take $ \mu = dx \wedge d y \wedge dz $ to be the standard volume of $ \mathbb{R}  ^3 $ restricted to $ \mathcal{D} $.
We shall denote by $ \operatorname{div}_ \mu ( \mathbf{u} )$ the divergence of $ \mathbf{u} $ with respect to the volume form $ \mu $. 
Recall that the divergence is the function $\operatorname{div}_ \mu ( \mathbf{u} )$ defined by the equality
\[
\pounds _ \mathbf{u} \mu = \operatorname{div}_ \mu ( \mathbf{u} ) \mu , 
\]
in which $ \pounds _ \mathbf{u} $ denotes the Lie derivative with respect to the vector field 
$\mathbf{u} $, see, e.g., \cite{AbMa1978}. When $\mu $ is the standard volume, one evidently recovers 
the usual divergence operator $\operatorname{div}$ on vector fields.

\paragraph{Diffeomorphism groups.}
Let us denote by $ \operatorname{Diff}_{ \mu }( \mathcal{D} )$ the group of all smooth diffeomorphisms 
$ \varphi : \mathcal{D} \rightarrow \mathcal{D} $ that preserve the volume form $ \mu $, i.e., $ \varphi^\ast \mu =\mu $. 
The group structure on $ \operatorname{Diff}_{ \mu }( \mathcal{D} )$ is given by the composition of diffeomorphisms. 
The group $ \operatorname{Diff}_{ \mu }( \mathcal{D} )$ can be endowed with the structure of a Fr\'echet Lie group, 
although in this paper we shall only use the Lie group structure at a formal level. The Lie algebra of the group 
$ \operatorname{Diff}_{ \mu }( \mathcal{D} )$ is given by the space $ \mathfrak{X} _ \mu ( \mathcal{D} )$ of all 
divergence free (relative to $\mu$) vector fields on $ \mathcal{D} $, parallel to the boundary $ \partial \mathcal{D} $:
\[
\mathfrak{X}  _ \mu ( \mathcal{D} )= \{ \mathbf{u} \in \mathfrak{X}  ( \mathcal{D} ) \mid \operatorname{div}_ \mu ( \mathbf{u} )=0, \;\; \mathbf{u}\, \|\, \partial \mathcal{D} \}.
\]

Given the strictly positive function $\bar \sigma >0$ on $ \mathcal{D} $, we consider the new volume form $ \bar \sigma \mu $ 
with associated diffeomorphism group and Lie algebra denoted $  \operatorname{Diff}_{ \bar \sigma \mu }( \mathcal{D} )$ and 
$\mathfrak{X}  _{ \bar \sigma \mu} ( \mathcal{D} )= \{ \mathbf{u} \in \mathfrak{X}  ( \mathcal{D} ) \mid 
\operatorname{div}_ {  \bar \sigma \mu}( \mathbf{u} )=0, \;\; \mathbf{u}\, \|\, \partial \mathcal{D} \}$, respectively.
In the next Lemma, we rewrite the condition $\operatorname{div}_ {  \bar \sigma \mu}( \mathbf{u} )=0$ by using exclusively 
the divergence operator $ \operatorname{div}_ \mu $ associated to the initial volume form $ \mu $.

\begin{lemma}\label{Lemma1}  Let $ \mathcal{D}$ be a manifold endowed with a volume form $\mu $ and let $\bar \sigma >0$ be a 
strictly positive smooth function on $ \mathcal{D} $. Then we have
\[
\operatorname{div}_ \mu (\bar \sigma  \mathbf{u} )= \bar \sigma  \operatorname{div}_{\bar \sigma  \mu } (\mathbf{u} ).  
\]
\end{lemma} 
\begin{proof}
 We will use the following properties of the Lie derivative $ \pounds _ \mathbf{u} $, 
  the exterior differential $ \mathbf{d} $, and the inner product $ \mathbf{i}_ \mathbf{u} $ on differential 
  forms (see, e.g., \cite{AbMa1978}): for a $k$-form $ \alpha $, an $n$-form $ \beta $, and a vector field $ \mathbf{u} $, 
  we have 
  \begin{align*}\notag 
 &  \pounds _\mathbf{u} \alpha = \mathbf{d} \left( \mathbf{i} _ \mathbf{u} \alpha \right) + \mathbf{i} _ \mathbf{u} \mathbf{d} \alpha , 
  \quad \mathbf{d} ( \alpha \wedge \beta )= \mathbf{d} \alpha \wedge \beta  +(-1) ^k  \alpha \wedge \mathbf{d} \beta, \\
 & \mathbf{i} _ \mathbf{u} ( \alpha \wedge \beta ) = \mathbf{i} _ \mathbf{u} \alpha \wedge \beta +(-1) ^k  \alpha \wedge \mathbf{i} _ \mathbf{u} \beta.
  \end{align*}
  On the one hand, we have
  \begin{align*} 
  \operatorname{div}_ \mu (\bar \sigma \mathbf{u} )\mu &= \pounds _{ \bar \sigma \mathbf{u} } \mu = \mathbf{d} \left( \mathbf{i} _{ \bar \sigma \mathbf{u} } \mu \right)= 
  \mathbf{d} \left( \bar \sigma \mathbf{i} _ \mathbf{u} \mu \right) =\mathbf{d} \bar \sigma \wedge  \mathbf{i} _ \mathbf{u} \mu 
  + \bar \sigma \mathbf{d} \left( \mathbf{i} _ \mathbf{u} \mu \right)\\
  &=\left(  \mathbf{i} _ \mathbf{u} \mathbf{d} \bar \sigma\right) \mu - \mathbf{i} _ \mathbf{u} \left( \mathbf{d} \bar \sigma \wedge \mu \right) 
  + \bar \sigma \operatorname{div}_ \mu \mathbf{u}  = \left( \mathbf{d} \bar \sigma \cdot \mathbf{u} \right) \mu +\bar \sigma \operatorname{div}_ \mu \mathbf{u}. 
  \end{align*} 
  On the other hand, we have
  \[
  \bar \sigma\operatorname{div}_{\bar \sigma \mu } (\mathbf{u} ) \mu = \pounds _ \mathbf{u} ( \bar \sigma \mu )=\left( \mathbf{d} \bar \sigma \cdot \mathbf{u} \right) \mu 
  +\bar \sigma \pounds _ \mathbf{u} \mu = \left( \mathbf{d} \bar \sigma \cdot \mathbf{u} \right) \mu + \bar \sigma \operatorname{div}_ \mu ( \mathbf{u} ).
  \]
  This proves the result.
\end{proof}

\medskip

From this Lemma, we deduce that the appropriate Lie groups associated to the anelastic and pseudo-incompressible systems are given by
\[
G= \operatorname{Diff}_{ \bar\rho \mu }( \mathcal{D} ) \quad\text{and}\quad  G= \operatorname{Diff}_{\bar \rho  \bar\theta \mu }( \mathcal{D} ),
\]
respectively. Indeed, from the preceding Lemma, it follows that the Lie algebras of these groups can be written as 
\begin{align*} 
\mathfrak{X}  _{ \bar \rho  \mu} ( \mathcal{D} )&= \{ \mathbf{u} \in \mathfrak{X}  ( \mathcal{D} ) \mid \operatorname{div}_ {  \mu}( \bar \rho \mathbf{u} )=0, \;\; 
\mathbf{u}\, \|\, \partial \mathcal{D} \} \;\;\; \text{and}\\
\mathfrak{X}  _{ \bar\rho \bar\theta  \mu} ( \mathcal{D} )&= \{ \mathbf{u} \in \mathfrak{X}  ( \mathcal{D} ) \mid \operatorname{div}_ {  \mu}( \bar\rho \bar\theta \mathbf{u} )=0, \;\; \mathbf{u}\, \|\, \partial \mathcal{D} \},
\end{align*} 
respectively. They correspond to the anelastic and pseudo-incompressible constraints on the fluid velocity.
We will continue to uses the subscript $\bar \sigma\mu$ when referring to both Lie groups and both Lie algebras.

\paragraph{Euler-Poincar\'e variational principles.} The diffeomorphism group $ \operatorname{Diff}_{\bar \sigma \mu }( \mathcal{D} )$ plays the 
role of the configuration manifold for these fluid models. The motion of the fluid is completely characterized by a time 
dependent curve $\varphi (t,\_\,)\in \operatorname{Diff}_{ \bar \sigma  \mu }( \mathcal{D} )$: a particle located at a point $X \in \mathcal{D}  $ 
at time $t=0$ travels to ${\bf x} = \varphi (t,X) \in \mathcal{D} $ at time $t$.
Exactly as in classical mechanics, the Lagrangian of the system is defined on the tangent bundle $T  \operatorname{Diff}_{ \bar \sigma  \mu }( \mathcal{D} )$ 
of the configuration manifold. We shall denote it by $L_{ \Theta _0}: T  \operatorname{Diff}_{ \bar \sigma  \mu }( \mathcal{D} ) \rightarrow \mathbb{R} $. 
The index $ \Theta _0$ indicates that this Lagrangian parametrically depends on the potential temperature $ \Theta _0(X)$ 
that is expressed here in the Lagrangian description.

The equations of motion in the Lagrangian description follow from the Hamilton principle
\begin{equation}\label{HP} 
\delta \int_0^T L_{ \Theta  _0} ( \varphi , \dot \varphi )dt=0,
\end{equation}
over a time interval $[0,T]$, for variations $\delta \varphi$ with $ \delta\varphi (0)= \delta \varphi (T)=0$.

In the Eulerian description, the variables are the Eulerian velocity $ \mathbf{u}(t, \mathbf{x} ) $ and 
the potential temperature $ \theta (t, \mathbf{x} )$. 
They are related to $\varphi (t,X)$ and $ \Theta _0(X)$ as
\begin{equation}\label{def_Eulerian} 
\mathbf{u} (t, \varphi (t,X))= \dot\varphi (t,X) \quad\text{and}\quad \theta (t, \varphi (t,X))=\Theta _0(X).
\end{equation} 
We assume that the Lagrangian $L_{ \Theta _0}$ can be rewritten exclusively in terms of these two Eulerian 
variables, and we denote it by $\ell( \mathbf{u} , \theta )$. This assumption means that $L_{ \Theta  _0}$ is right-invariant with respect to the action of the subgroup
\[
\operatorname{Diff}_{ \bar \sigma  \mu }( \mathcal{D} )_{ \Theta _0}=\{ \varphi \in \operatorname{Diff}_{ \bar \sigma  \mu }( \mathcal{D} )\mid \Theta _0 (\varphi (X))=\Theta _0(X),\;\forall\, X \in \mathcal{D} \}
\]
of all diffeomorphisms that keep $ \Theta _0 $ invariant.

By rewriting the Hamilton principle \eqref{HP} in terms of the Eulerian variables $ \mathbf{u} $ and $ \theta $, we get the Euler-Poincar\'e variational principle 
\begin{equation}\label{EP_VP} 
\delta \int_0^T\ell( \mathbf{u} , \theta )dt=0, \quad\text{for variations} \quad\delta  \mathbf{u} 
= \partial _t \mathbf{v} +[ \mathbf{u} ,\mathbf{v} ], \quad \delta \theta =- \mathbf{d}  \theta \cdot \mathbf{v},
\end{equation} 
where $ \mathbf{v} (t, \mathbf{x} )$ is an arbitrary vector field on $ \mathcal{D} $ parallel to the boundary and with $ \operatorname{div} _\mu (\bar \sigma\mathbf{v} )=0$, (i.e., $ \mathbf{v} \in \mathfrak{X}  _{\bar \sigma \mu }( \mathcal{D} )$ by Lemma \ref{Lemma1}), and with $ \mathbf{v} (0, \mathbf{x} )= \mathbf{v} (T, \mathbf{x} )=0$. The bracket $[ \mathbf{u} , \mathbf{v} ]$, locally given by $[ \mathbf{u} ,\mathbf{v} ]^i:= u ^j \partial _j v ^i - v ^j \partial _j  u ^i$, is the Lie bracket of vector fields.

The expressions for $ \delta \mathbf{u} $ and $ \delta \theta $ in \eqref{EP_VP} follow by taking the variation of the first and second 
equalities in \eqref{def_Eulerian} and defining $ \mathbf{v} (t, \mathbf{x} )$ as $ \mathbf{v} (t,\varphi (t,X))= \delta\varphi (t,X)$. A direct and efficient way to obtain these expressions, or the variational principle \eqref{EP_VP},
is to apply the general theory of Euler-Poincar\'e reduction on Lie groups, see \cite{HoMaRa1998}.

In order to compute the associated equations, one needs to fix an appropriate space in nondegenerate duality with $\mathfrak{X}  _{\bar \sigma \mu }( \mathcal{D} )$. 
This is recalled in the next Lemma, which follows from the Hodge decomposition and shall play a crucial role in the discrete setting later. 
Recall that given a vector space $V$, a space in nondegenerate duality with $V$ is a vector space $V'$  
together with a bilinear form $ \left\langle \;\,, \;\right\rangle : V' \times V \rightarrow \mathbb{R} $ such that 
$ \left\langle \alpha , v \right\rangle =0$, for all $v \in V$, implies $ \alpha =0$ and $ \left\langle \alpha , v \right\rangle =0$, 
for all $ \alpha \in V'$, implies $v=0$.

\begin{lemma}\label{Lemma2} The space $ \Omega ^1 (\mathcal{D} )/ \mathbf{d} \Omega ^0 ( \mathcal{D} )$ of one-forms modulo exact forms is in 
nondegenerate duality with the space $ \mathfrak{X}_{ \bar\sigma \mu }( \mathcal{D} )$, the Lie algebra of $ \operatorname{Diff}_{\bar \sigma \mu }( \mathcal{D} )$. 
The nondegenerate duality pairing is given by
\begin{equation}\label{pairing_vector} 
\left\langle \;\,,\; \right\rangle _{\bar \sigma } :\Omega ^1 (\mathcal{D} )/ \mathbf{d} \Omega ^0 ( \mathcal{D} ) \times  \mathfrak{X}_{\bar \sigma \mu }( \mathcal{D} ) \rightarrow \mathbb{R}  , 
 \quad\left\langle [ \alpha ], \mathbf{u}  \right\rangle _{\bar \sigma} := \int_\mathcal{D} (\alpha \cdot\mathbf{u}) \bar \sigma\mu,
\end{equation} 
where $[ \alpha ]$ denotes the equivalence class of $ \alpha $ modulo exact forms.
\end{lemma}
 \begin{proof}
 It is well-known that if $g$ is a Riemannian metric, with $\mu _g$ the associated volume form on $ \mathcal{D} $, then
  \[
  \left\langle \;\,, \; \right\rangle : \Omega ^1 ( \mathcal{D} )/ \mathbf{d} \Omega ^0 ( \mathcal{D} ) \times \mathfrak{X}  _{ \mu _g}( \mathcal{D} ) \rightarrow \mathbb{R}  , 
  \quad \left\langle [ \alpha ], \mathbf{v} \right\rangle =\int_\mathcal{D}(\alpha \cdot \mathbf{v} )\mu_g , 
  \]
  is a nondegenerate duality pairing, see e.g., \cite[\S14.1]{MaRa1994}. This result follows from the Hodge decomposition of $1$-forms, which needs the introduction of a Riemannian metric $g$.

  In our case, the volume forms $ \mu $ and $ \bar \sigma \mu $ are not necessarily associated to a Riemannian metric. We shall thus introduce a Riemannian metric $g$ uniquely for the purpose of this proof, with associated Riemannian volume form $ \mu _g$. Let $ f$ be the function defined by $ \bar\sigma  \mu = f\mu _g$. Since $ \mathcal{D} $ is orientable and connected, we have either $f>0$ or $ f<0$ on $ \mathcal{D} $. We can rewrite the duality pairing \eqref{pairing_vector} as
  \begin{equation}\label{ND_DP} 
  \int _ \mathcal{D} (\alpha \cdot\mathbf{u}) \bar \sigma\mu= \int_ \mathcal{D} ( \alpha \cdot (f \mathbf{u} )) \mu _g.
  \end{equation} 
  By successive applications of Lemma \ref{Lemma1}, we have
  \[
  \operatorname{div}_{ \mu _g }( f \mathbf{u} )=\operatorname{div}_{\frac{ \bar\sigma }{f} \mu } ( f\mathbf{u} )= f\operatorname{div}_{\bar\sigma \mu }(\mathbf{u} )= \frac{f}{\bar \sigma }\operatorname{div}_ \mu ( \bar \sigma \mathbf{u} )=0,      
  \]
  where the last equality follows since $ \mathbf{u} \in \mathfrak{X}  _{ \bar \sigma \mu }( \mathcal{D} )$. 
  This proves that $\mathbf{v}= f \mathbf{u} \in\mathfrak{X}  _{ \mu _g}( \mathcal{D} )$. 
  We can thus write the duality pairing $ \left\langle \;\,,\; \right\rangle _{\bar \sigma }$ 
  in terms of the nondegenerate duality pairing \eqref{ND_DP} as $ \left\langle [ \alpha ],  \mathbf{u}\right\rangle _{\bar \sigma } = \left\langle [ \alpha ], f \mathbf{u} \right\rangle $, 
  which proves that it is nondegenerate. 
\end{proof}

\bigskip

In a similar way to  \eqref{pairing_vector}, we shall identify the dual to the space of functions $ \mathcal{F} (\mathcal{D} )$ 
with itself by using the nondegenerate duality pairing
\begin{equation}\label{pairing_function}
\mathcal{F}(\mathcal{D} ) \times \mathcal{F} ( \mathcal{D} ) \rightarrow \mathbb{R}  , 
\quad \left\langle h, \theta  \right\rangle_{\bar \sigma}  =\int_ \mathcal{D} (h \theta) \bar \sigma \mu .
\end{equation}

Given a Lagrangian $\ell: \mathfrak{X}  _{ \bar \sigma \mu}( \mathcal{D} ) \times \mathcal{F} ( \mathcal{D} ) \rightarrow \mathbb{R}$, 
the functional derivatives of $\ell$ are defined with respect to the parings \eqref{pairing_vector} and \eqref{pairing_function} and denoted
\[
\left[ \frac{\delta \ell}{\delta \mathbf{u} } \right]  \in \Omega ^1 (\mathcal{D} )/ \mathbf{d} \Omega ^0 ( \mathcal{D} ), \;\; 
\text{for $ \frac{\delta \ell}{\delta \mathbf{u} } \in \Omega ^1 ( \mathcal{D} )$}, \qquad\text{and}\qquad \frac{\delta \ell}{\delta \theta } \in \mathcal{F} ( \mathcal{D} ).
\]

\begin{proposition}\label{prop_EP}  The variational principle \eqref{EP_VP} yields the partial differential equation
\begin{equation}\label{EP_general} 
\partial _t \frac{\delta \ell}{\delta \mathbf{u} }+\pounds _\mathbf{u} \frac{\delta \ell}{\delta \mathbf{u} } + \frac{\delta \ell}{\delta \theta } \mathbf{d} \theta = - \mathbf{d} p, \quad \text{with} \quad \operatorname{div}_ \mu (\bar \sigma\mathbf{u} )=0, \;\; \mathbf{u}\, \|\, \partial \mathcal{D},
\end{equation}
where $ \pounds _ \mathbf{u} $ denotes the Lie derivative acting on one-forms, given by $ \pounds _ \mathbf{u} \alpha = \mathbf{d} (\mathbf{i} _ \mathbf{u} \alpha )+ \mathbf{i} _ \mathbf{u} \mathbf{d} \alpha $. This equation is supplemented with the advection equation
\[
\partial _t\theta + \mathbf{d}\theta\cdot \mathbf{u} =0,
\]
which follows from the definition of $ \theta $ in \eqref{def_Eulerian}. 
\end{proposition} 
\begin{proof}
 By definition of the functional derivatives, we have
  \[
  \delta \int_0^T\ell( \mathbf{u} , \theta ) dt= \int_0^T\int_ \mathcal{D} \frac{\delta \ell}{\delta\mathbf{u} }\cdot \delta \mathbf{u}\,  \bar \sigma \mu dt+ \int_0^T\int_ \mathcal{D} \frac{\delta \ell}{\delta \theta  }\cdot \delta \theta  \, \bar \sigma \mu dt.
  \]
  Using the expression for $ \delta \mathbf{u} $ in \eqref{EP_VP}, integrating by parts, and using the equalities $ \pounds _ \mathbf{u} \mathbf{v} =[\mathbf{u} , \mathbf{v} ]$ and $ \mathbf{d} ( \alpha \cdot \mathbf{v} )\cdot \mathbf{u} =(\pounds_\mathbf{u} \alpha)\cdot \mathbf{v} + \alpha \cdot (\pounds _\mathbf{u} \mathbf{v} )$, the first term reads
  \[
  -\int_0^T \int_ \mathcal{D} \left(  \partial _t \frac{\delta \ell}{\delta\mathbf{u} } + \pounds _\mathbf{u}  \frac{\delta \ell}{\delta\mathbf{u} } \right)\cdot \mathbf{v} \bar \sigma \mu dt + \int_0^T \int_ \mathcal{D} \mathbf{d} \left( \frac{\delta \ell}{\delta \mathbf{u} } \cdot \mathbf{v} \right)\cdot \mathbf{u} \bar \sigma \mu dt.
  \]
  We can write $ \mathbf{d} \left( \frac{\delta \ell}{\delta \mathbf{u} } \cdot \mathbf{v} \right)\cdot \mathbf{u} 
    = \operatorname{div}_{\bar \sigma \mu }  \left(  \frac{\delta \ell}{\delta \mathbf{u} } \cdot \mathbf{v} \, \mathbf{u} \right) 
    -  \frac{\delta \ell}{\delta \mathbf{u} } \cdot \mathbf{v} \operatorname{div}_{\bar \sigma \mu }(\mathbf{u} )
    =\operatorname{div}_{\bar \sigma \mu }  \left(  \frac{\delta \ell}{\delta \mathbf{u} } \cdot \mathbf{v} \, \mathbf{u} \right)$, 
    since $ \operatorname{div}_{\bar \sigma \mu }(\mathbf{u} ) =0$. Then, by the Gauss Theorem,
  \[
  \int_ \mathcal{D} \operatorname{div}_{\bar \sigma \mu }  \left(  \frac{\delta \ell}{\delta \mathbf{u} } \cdot \mathbf{v} \, \mathbf{u} \right)\bar \sigma \mu
  =\int_{\partial \mathcal{D} }\frac{\delta \ell}{\delta \mathbf{u} } \cdot \mathbf{v}  \,   \mathbf{i} _ \mathbf{u}\mu\, \bar \sigma=0 ,
  \]
  since $ \mathbf{u} \,\|\, \partial\mathcal{D} $. Combining these results, we thus get
  \[
  \int_0^T \int_\mathcal{D} \left( \partial _t \frac{\delta \ell}{\delta \mathbf{u} }+\pounds _\mathbf{u} \frac{\delta \ell}{\delta \mathbf{u} } + \frac{\delta \ell}{\delta \theta } \mathbf{d} \theta  \right) \cdot \mathbf{v} \bar \sigma \mu dt=0,
  \]
  for all $ \mathbf{v} \in \mathfrak{X}  _{\bar \sigma \mu }( \mathcal{D} )$. By Lemma \ref{Lemma2}, 
  it follows that the one-form $\partial _t \frac{\delta \ell}{\delta \mathbf{u} }+\pounds _\mathbf{u} \frac{\delta \ell}{\delta \mathbf{u} } + \frac{\delta \ell}{\delta \theta } \mathbf{d} \theta  $ 
  is exact, i.e., there exists a function $p$ such that this expression equals $ \mathbf{d} p$.
  \end{proof}

\medskip

\begin{remark}{\rm We note that the statement of Proposition \ref{prop_EP} does not need the introduction of a Riemannian metric 
$g$ on $ \mathcal{D} $. Only a volume form $ \mu $ is fixed, together with a strictly positive function $\bar \sigma $. It can be 
however advantageous to formulate the equations \eqref{EP_general} in terms of a Riemannian metric $g$ (note that we do not suppose 
that $ \mu $ or $\bar \sigma \mu $ equals $\mu _g$). In this case, identifying one-forms and vector fields via the flat operator 
$\mathbf{u} \in  \mathfrak{X}  ( \mathcal{D} ) \rightarrow \mathbf{u} ^\flat = g( \mathbf{u} , \_\,) \in \Omega ^1 (\mathcal{D} )$, 
the space $ \Omega ^1 (\mathcal{D} )/ \mathbf{d} \Omega ^0 ( \mathcal{D} )$ can be identified with the space of vector fields $\mathfrak{X}  ( \mathcal{D} )$ modulo gradient (with respect to $g$) of functions. The nondegenerate duality pairing \eqref{pairing_vector} thus reads
\begin{equation}\label{pairing_g} 
\left\langle [ \mathbf{v} ], \mathbf{u} \right\rangle _{\bar \sigma }= \int_ \mathcal{D} g( \mathbf{v} , \mathbf{u} )\bar \sigma \mu .
\end{equation} 
In terms of this duality pairing, the equations \eqref{EP_general} are equivalently written as
\begin{equation}\label{EP_general2}
\partial _t \frac{\delta \ell}{\delta \mathbf{u} } + \mathbf{u}  \cdot \nabla \frac{\delta \ell}{\delta \mathbf{u} }
+\nabla\mathbf{u}  ^\mathsf{T} \cdot   \frac{\delta \ell}{\delta \mathbf{u} } + \frac{\delta \ell}{\delta \theta }\nabla  \theta=- \nabla p, 
\end{equation}
where $ \nabla $ acting on a vector field is the covariant derivative associated to the Riemannian metric $g$, $\nabla $ acting on a function is the gradient relative to $g$, and $ \nabla \mathbf{u} ^\mathsf{T}$ denotes the transpose with respect to $g$.
}
\end{remark}

% 
% \todo{\incl{WB: It is unclear, in which setting we are now for the following derivations: 
% are we now in $\mathbbm R^3$ with Riemannian metric. but if we use $\mu$ instead of $\bar \sigma \mu$ 
% how can we derive the anelastic and pseudo-incompressible equations?}\\
% \color{magenta}FGB: We are in a domain in $ \mathbb{R}  ^2 $ or a domain in $ \mathbb{R}  ^3 $ with smooth 
%boundary. Except this, we are in the same setting as before: we fix a volume form $ \mu $, and we work on the 
%space $ \mathfrak{X}_{\bar \sigma \mu }( \mathcal{D} )$. The only difference with before is that now there is a  
%Riemannian metric, chosen to be the inner product on $ \mathbb{R}  ^2 $ or on $ \mathbb{R}  ^3 $, which appears 
%in the kinetic energy of the Lagrangian. Note that $\mu $ is not necessarily related to $g$.\\
% I added comments about this below.
% }

We shall now apply this setting to the anelastic and the pseudo-incompressible equations. The fluid domain $ \mathcal{D} $ is a subset of the vertical plane $ \mathbb{R}  ^2 \ni (x,z)$ or of the space $ \mathbb{R}  ^3 \ni (x,y,z)$, and has a smooth boundary $ \partial\mathcal{D} $. We fix a volume form $ \mu $ on $ \mathcal{D} $.

\paragraph{1) Anelastic equations.} For the anelastic equation, we take $ \bar \sigma (z)= \bar \rho (z)$, the reference mass density. The Lagrangian is given by
\begin{equation}\label{Lagrangian_anelastic} 
\ell( \mathbf{u} , \theta )= \int_{ \mathcal{D} }\left(  \frac{1}{2} | \mathbf{u} | ^2 - c _p \bar \pi \theta \right) \bar \rho  \mu , 
\quad\mathbf{u} \in \mathfrak{X} _{\bar \rho \mu }( \mathcal{D} ),
\end{equation}
where $ \bar \pi (z)$ is such that $c _p \frac{\partial \bar \pi }{\partial z}=- \frac{g}{\bar \theta}$ and the norm is computed relative to the standard inner product on $ \mathbb{R}  ^2 $ or $ \mathbb{R}  ^3 $.

Relative to the pairings \eqref{pairing_g} and \eqref{pairing_function} we get
\begin{equation}\label{functional_derivative_A} 
\frac{\delta \ell}{\delta \mathbf{u} }= \mathbf{u} \quad \text{and} \quad \frac{\delta \ell}{\delta \theta}=  - c _p \bar \pi,
\end{equation} 
so that the Euler-Poincar\'e equations \eqref{EP_general2} read 
$\partial _t \mathbf{u} + \mathbf{u} \cdot \nabla \mathbf{u} +  \nabla \mathbf{u} ^\mathsf{T} \cdot \mathbf{u}-  c_ p \bar \pi \nabla \theta  = - \nabla p$, in terms of the pressure $p$.
To permit a comparison of these anelastic equations with those given in the standard form of \eqref{anelastic}
in terms of Exner pressure $\pi'$, we note that $ \nabla \mathbf{u} ^\mathsf{T} \cdot \mathbf{u} = \frac{1}{2} \nabla | \mathbf{u} | ^2 $ and that 
$-c _p \bar \pi \nabla \theta $ differs from  $-  g\frac{\theta' }{\bar \theta }\mathbf{z}$ by a gradient term, indeed:
\[
-c _p \bar \pi \nabla \theta  = - c _p \nabla \left( \bar \pi \theta  \right) + c _p (\nabla \bar \pi) \theta = - c _p \nabla \left( \bar \pi \theta  \right) -  g\frac{\theta }{\bar \theta }\mathbf{z}  = - \nabla \left( c _p \bar \pi \theta + gz \right)-  g\frac{\theta' }{\bar \theta }\mathbf{z}.
\]
Therefore, with $ \pi '$ defined in terms of $p$ by the equality $c_p \bar \theta \pi '=p+ \frac{1}{2} | \mathbf{u} |^2 -g\mathbf{z} -c_p \bar \pi \theta $, the Euler-Poincar\'e equations yield the anelastic equations \eqref{anelastic}. 

\paragraph{2) Pseudo-incompressible equations.} In this case we take $ \bar \sigma (z) := \bar \rho (z)  \bar \theta (z)$ and the Lagrangian is given by
\begin{equation}\label{Lagrangian_pseudo_incomp} 
\ell( \mathbf{u} , \theta )= \int_{ \mathcal{D} } \frac{1}{ \theta }\left(  \frac{1}{2} | \mathbf{u} | ^2 - gz \right) \bar \rho  \bar \theta  \mu ,
\quad\mathbf{u} \in \mathfrak{X} _{\bar \rho \bar \theta \mu }( \mathcal{D} ).
\end{equation} 
As before, the kinetic energy is computed relative to the standard inner product on $ \mathbb{R}  ^2 $ or $ \mathbb{R}  ^3 $. Relative to the pairings \eqref{pairing_g} and \eqref{pairing_function} we get
\begin{equation}\label{functional_derivative_PI} 
\frac{\delta \ell}{\delta \mathbf{u} }= \frac{1}{ \theta } \mathbf{u} \quad \text{and} \quad \frac{\delta \ell}{\delta \theta}=  - \frac{1}{\theta ^2} \left( \frac{1}{2} | \mathbf{u} | ^2 -gz \right) ,
\end{equation}
so that the Euler-Poincar\'e equations \eqref{EP_general2} read
\[
\partial _t \left( \frac{1}{\theta } \mathbf{u}\right)  + \mathbf{u} \cdot \nabla \left( \frac{1}{\theta }\mathbf{u}\right) +\nabla \mathbf{u}^\mathsf{T} \cdot  \frac{1}{\theta } \mathbf{u}-\frac{1}{\theta ^2 }\left(  \frac{1}{2} | \mathbf{u} | ^2 - g z  \right)\nabla \theta =- \nabla p.
\]
After some computations, using the relation $- \frac{g}{\bar \theta }=c_p \partial _z \bar \pi $,  these equations recover the pseudo-incompressible system \eqref{pseudo_incompressible} with $c_p( \bar\pi + \pi ')= p + \frac{1}{\theta } ( \frac{1}{2} | \mathbf{u} | ^2 -g z)$.
% \todo{\incl{Derivation not completely clear to me.}\\
% \color{magenta}FGB: Ok, this is a little bit long, but it only uses vector calculus so it cannot be written. We shall discuss it.}

\medskip

Based on these results, we can formulate the following statement that will allow us to derive the variational discretization 
of these two models by the discrete diffeomorphism group approach.

\begin{theorem}\label{crucial}  Consider a domain $ \mathcal{D} $ with smooth boundary $ \partial \mathcal{D}$ and volume form $ \mu $. The anelastic system with reference density $\bar \rho $, resp., the pseudo-incompressible system with reference density $\bar \rho $ and reference potential temperature $\bar \theta $ can be derived from an Euler-Poincar\'e variational principle for the Lie group
\begin{equation}\label{Lie_groups} 
G= \operatorname{Diff}_{ \bar\rho  \mu }( \mathcal{D} )\quad\text{resp.}\quad   G= \operatorname{Diff}_{ \bar\rho  \bar\theta \mu }( \mathcal{D} ),
\end{equation} 
with Lagrangian \eqref{Lagrangian_anelastic}, resp., \eqref{Lagrangian_pseudo_incomp}.
\end{theorem} 

\paragraph{Kelvin-Noether circulation theorems.} The Euler-Poincar\'e formulation is well adapted for a systematic derivation of the circulation theorems, see \cite{HoMaRa1998}. 
From \eqref{EP_general} one indeed deduces the following general form of the circulation theorem
\begin{equation}\label{KN1} 
\frac{d}{dt} \oint_{c_t} \frac{\delta \ell}{\delta \mathbf{u} }  =-  \oint_{c_t}  \frac{\delta \ell}{\delta \theta } \mathbf{d}   \theta ,
\end{equation} 
where $c_t= \varphi (t, c_0)$ is a loop advected by the fluid flow $ \varphi (t,\_\,) $ and $\oint_{ c_t} \alpha $ denotes the circulation of the one-form $ \alpha $ along the loop $c_t$.
Using the equation $ \partial _t \theta + \mathbf{d} \theta  \cdot \mathbf{u}   =0$, one also deduces another useful form, namely,
\begin{equation}\label{KN2} 
\frac{d}{dt} \oint_{c_t} \theta \frac{\delta \ell}{\delta \mathbf{u} } =-  \oint_{c_t} \left(  \theta  \frac{\delta \ell}{\delta \theta } \mathbf{d}  \theta- \theta \mathbf{d}  p \right) .
\end{equation} 
We shall not present the derivation of \eqref{KN1} and \eqref{KN2} since they follow from similar arguments with those explained in details in \cite{HoMaRa1998}.

For the anelastic system, using \eqref{functional_derivative_A} and the equalities $ c_p \bar \pi \mathbf{d}  \theta = c_p \mathbf{d}  ( \bar \pi \theta )- c_p \mathbf{d}   \bar \pi \theta = c_p \mathbf{d}  ( \bar \pi \theta )+ g \frac{ \theta }{\bar \theta } \mathbf{z} $, expression \eqref{KN1} yields the equivalent forms
\[
\frac{d}{dt} \oint_{c_t}  \mathbf{u} \cdot d \mathbf{x}   = c_p \oint_{c_t}\bar \pi \mathbf{d} \theta   \qquad\text{or}\qquad \frac{d}{dt} \oint_{c_t}  \mathbf{u} \cdot d \mathbf{x}  = g \oint_{c_t} \frac{\theta }{\bar \theta } \mathbf{z} \cdot d \mathbf{x} 
\] 
of the circulation theorem.
For the pseudo-incompressible system, using \eqref{functional_derivative_PI} and the equalities $ \frac{1}{ \theta } \left( \frac{1}{2} | \mathbf{u} | ^2 -g z \right) \mathbf{d}  \theta -\theta \mathbf{d}  p= \mathbf{d}  \left( \frac{1}{2} | \mathbf{u} | ^2 -g z \right) - c_p \mathbf{d}  \pi\theta  $, the expression \eqref{KN2} yields
\[
\frac{d}{dt} \oint_{c_t}  \mathbf{u} \cdot d \mathbf{x} = c_p \oint_{c_t} \pi \mathbf{d} \theta .
\]
As shown further below, these conservation laws of the continuous equations, here presented with explicit formulas, 
are also preserved by the discrete variational discretizations that we will derive in the following section.

\section{Variational discretizations}\label{VD}

In this section we first quickly review from \cite{PaMuToKaMaDe2010} the discrete diffeomorphism group approach 
in the incompressible case. Then, based on the results of Theorem \ref{crucial}, we show that an appropriate 
adaptation of this approach allows us to derive a variational discretization of the anelastic and pseudo-incompressible 
systems valid on a large class of mesh discretizations of the fluid domain.

\paragraph{Review of the discrete diffeomorphism group approach in the incompressible case.}
The spatial discretization of the equations is accomplished by considering the finite dimensional approximation of the 
group of volume preserving diffeomorphisms developed in \cite{PaMuToKaMaDe2010}, which we roughly recall below.

Given a mesh $\mathbb{M}$ on the fluid domain $ \mathcal{D} $ with cells $C_i$, $i=1,...,N$, define a diagonal $N\times N$ matrix $\Omega$ consisting of cell 
volumes: $\Omega_i=\operatorname{Vol}(C_i)$. The discretization of the group $ \operatorname{Diff}_{ \mu }( \mathcal{D} )$ of volume preserving 
diffeomorphisms of $ \mathcal{D} $ is the matrix group
\begin{equation}\label{DD}
\mathsf{D}(\mathbb{M})=\left\{q\in \operatorname{GL}(N)^+\mid q\cdot\mathbf{1}=\mathbf{1}\;\;\text{and}\;\; q^\mathsf{T}\Omega q=\Omega\right\},
\end{equation}
where $\operatorname{GL}(N)^+$ is the group of invertible $N \times N$ matrices with positive determinant, and $\mathbf{1}$ denotes the column $(1,...,1)^\mathsf{T}$ so that the first condition reads $\sum_{j=1}^Nq_{ij}=1$ for all $i=1,...,N$. 
The main idea behind this definition is the following (see \cite{PaMuToKaMaDe2010} for the detailed treatment). Consider the linear action of $\operatorname{Diff}_\mu ( \mathcal{D} )$ on the space $\mathcal{F}( \mathcal{D} )$ 
of functions on $ \mathcal{D} $, given by
\begin{equation}\label{action_diffeo}
f\in\mathcal{F}( \mathcal{D} )\mapsto f\circ\varphi^{-1}\in\mathcal{F}( \mathcal{D} ),\quad \varphi\in\operatorname{Diff}_\mu ( \mathcal{D} ).
\end{equation}
This linear map has two key proporties:\\
$(1)$ it preserves the $L^2$ inner product of functions;\\
$(2)$ it preserves constant functions $C$ on $ \mathcal{D} $: $C\circ\varphi^{-1}=C$.\\
In the discrete setting, a function is discretized as a vector $F\in\mathbb{R}^N$ whose value $F_i$ on cell $C_i$ is regarded as the cell average of the 
function. Accordingly, the discrete $L^2$ inner product of two discrete functions is defined by
\[
\left\langle F, G\right\rangle=F^\mathsf{T}\Omega G=\sum_{i=1}^NF_i\Omega_i G_i.
\]
The discrete diffeomorphism group \eqref{DD} is such that its action on discrete functions by matrix multiplication, is an approximation of the linear map \eqref{action_diffeo}. The conditions $q^\mathsf{T}\Omega q=\Omega$ and 
$q\cdot\mathbf{1}=\mathbf{1}$ are indeed the discrete analogues of the conditions $(1)$ and $(2)$ above, respectively.

The Lie algebra of $\mathsf{D}(\mathbb{M})$, denoted $\mathfrak{d}(\mathbb{M})$, is the space of $\Omega$-antisymmetric, row-null matrices:
\[
\mathfrak{d}(\mathbb{M})=\{A\in\mathfrak{gl}(N)\mid A\cdot \mathbf{1}=0\;\;\text{and}\;\; A^\mathsf{T}\Omega+\Omega A=0\}.
\]
The component $A _{ij} $ of the matrix $A$ is the weighted flux of the vector field $ \mathbf{u} $ through the face common 
to the cells $C_i$ and $C_j$. This relation induces a nonholonomic constraint on the Lie algebra $\mathfrak{d}(\mathbb{M})$ as only the fluxes through adjacent cells are non-zero:
\begin{equation}\label{nonholonomic_constraint}
\mathcal{S}  =\{A\in \mathfrak{d}(\mathbb{M})\mid A_{ij}\neq 0\Rightarrow j\in N(i)\},
\end{equation}
in which $N(i)$ is the set of all indices of cells adjacent to cell $C _i $.

Once a discrete Lagrangian $\ell_d: \mathfrak{d}(\mathbb{M}) \rightarrow \mathbb{R}$ has been selected, the derivation 
of the spatial variational discretization then proceeds by applying an Euler-Poincar\'e variational principle to this Lagrangian 
that takes into account the nonholonomic constraint. This approach has been developed in \cite{PaMuToKaMaDe2010} for the 
incompressible homogenous ideal fluid and extended to several models of incompressible fluids with advection equations 
in \cite{GaMuPaMaDe2010} and to rotating and/or stratified Boussinesq flows in \cite{DeGaGBZe2014}.

\paragraph{Discrete diffeomorphism groups for the two models.}
The results obtained in Theorem \ref{crucial} make it possible to extend this approach to treat the anelastic and pseudo-incompressible 
systems. The main idea consists in defining weighted versions of the volume of the cells in order to permit 
the use of the results recalled above in the incompressible case.

Given a mesh $ \mathbb{M}  $ on $ \mathcal{D} $, a reference density $ \bar \rho (z)$, and a reference potential temperature 
$\bar \theta (z)$ on $ \mathcal{D} $, we define, respectively, the diagonal matrices $ \Omega ^{\bar \rho }$ and $ \Omega ^{\bar \rho \bar \theta}$ of 
$ \bar \rho $\,-weighted and $ \bar \rho \bar \theta$\,-weighted volumes as
\[
\Omega ^{\bar \rho }_i := \int_{ C _i } \bar\rho (z)d \mathbf{x} \quad\text{and}\quad \Omega ^{\bar \rho \bar \theta}_i := \int_{ C _i } \bar\rho (z)\bar \theta(z)d \mathbf{x} .
\]
The discrete versions of the diffeomorphism groups $\operatorname{Diff}_{ \bar\rho  \mu }( \mathcal{D} )$ and $\operatorname{Diff}_{ \bar\rho  \bar\theta \mu }( \mathcal{D} )$ in \eqref{Lie_groups} are therefore 
\begin{align*} 
\mathsf{D}_{\bar \rho } ( \mathbb{M}  ):&= \big\{q\in \operatorname{GL}(N)^+\mid q\cdot\mathbf{1}=\mathbf{1}
\;\;\text{and}\;\; q^\mathsf{T}\Omega^{\bar \rho } q=\Omega^{\bar \rho }\big\},\\
\mathsf{D}_{\bar \rho\bar \theta } ( \mathbb{M}  ):& =\big\{q\in \operatorname{GL}(N)^+\mid q\cdot\mathbf{1}=\mathbf{1}
\;\;\text{and}\;\; q^\mathsf{T}\Omega ^{\bar \rho \bar \theta }q=\Omega^{\bar \rho \bar \theta }\big\},
\end{align*} 
with Lie algebras
\begin{align*} 
\mathfrak{d}_{\bar \rho } (\mathbb{M})&=\big \{A\in\mathfrak{gl}(N)\mid A\cdot \mathbf{1}=0\;\;\text{and}\;\; A^\mathsf{T}\Omega^{\bar \rho } +\Omega^{\bar \rho }  A=0\big \},\\
\mathfrak{d}_{\bar \rho \bar \theta } (\mathbb{M})&=\big \{A\in\mathfrak{gl}(N)\mid A\cdot \mathbf{1}=0\;\;\text{and}\;\; A^\mathsf{T}\Omega^{\bar \rho \bar \theta } +\Omega^{\bar \rho \bar \theta }  A=0\big \}.
\end{align*}
In order to treat the anelastic and pseudo-incompressible cases, we also need to appropriately modify the relation between the components $ A _{ij} $ and the velocity vector fields $ \mathbf{u} $, by taking into account the weights $ \bar \rho $ and $\bar \theta $. For $ \mathbf{u} \in \mathfrak{X}_{ \bar \rho \mu }( \mathcal{D})$, resp., 
$ \mathbf{u} \in \mathfrak{X}_{ \bar \rho   \bar \theta \mu }( \mathcal{D})$, we get
\begin{equation}\label{rel_A_u}
A_{ij}\simeq -\frac{1}{2\Omega^{\bar \rho }_i}\int_{D_{ij}}(\bar \rho \,\mathbf{u}\cdot\mathbf{n}_{ij})\mbox{d}S,\;\;\;\text{resp.},\;\;\;
A_{ij}\simeq -\frac{1}{2\Omega^{\bar \rho \bar\theta }_i}\int_{D_{ij}}(\bar\rho  \bar \theta \,\mathbf{u}\cdot\mathbf{n}_{ij})\mbox{d}S,
\end{equation}
in which $D_{ij}$ is the boundary common to $C_i$ and $C_j$, and $\mathbf{n}_{ij}$ is the normal vector field on $D_{ij}$ pointing from $C_i$ to $C_j$. The same nonholonomic constraint as before needs to be imposed, but this time on $ \mathfrak{d}_{\bar \sigma }( \mathbb{M}  )$, namely,
\begin{equation}\label{nonholonomic_constraint_space}
\mathcal{S}  _{ \bar\sigma }=\{A\in \mathfrak{d}_{\bar \sigma }(\mathbb{M})\mid A_{ij}\neq 0\Rightarrow j\in N(i)\}.
\end{equation}
This constraint induces a right-invariant linear constraint on the Lie group $\mathsf{D}_{\bar\sigma  } ( \mathbb{M}  )$: at $ q \in \mathsf{D}_{\bar \sigma  } (\mathbb{M})$, the constraint is defined as
\[
\mathcal{S}  _{ \bar\sigma }(q):=\{ \dot q \in T_q \mathsf{D}_{\bar \sigma  } (\mathbb{M}  )\mid \dot q q^{-1} \in \mathcal{S}  _{ \bar\sigma }\} \subset  T_q \mathsf{D}_{\bar \sigma  } (\mathbb{M}  ).
\]

In addition to the matrix $ A \in \mathfrak{d}_{\bar \sigma }( \mathbb{M}  )$ which discretizes the Eulerian velocity $ \mathbf{u} $, we introduce the discrete potential temperature $ \Theta \in\mathbb{R} ^N $ whose component $ \Theta _i$ is the average of the potential temperature $ \theta $ on cell $C _i $.

\paragraph{Variational discretization for the two models.}
The variational discretization is carried out by mimicking the Euler-Poincar\'e approach of Theorem \ref{crucial}.
Consider a discrete Lagrangian $L_{ \Theta _0,d}: T\mathsf{D}_{\bar \sigma  } (\mathbb{M}  ) \rightarrow \mathbb{R}  $ defined on the tangent bundle of the Lie group $\mathsf{D}_{\bar \sigma  } (\mathbb{M}  )$ and being an approximation of the Lagrangian in \eqref{HP}. The parameter $ \Theta _0 \in \mathbb{R}  ^N $ is the discrete potential temperature in the Lagrangian description. 

Hamilton's principle has to be appropriately modified to take into account the nonholonomic constraint, namely, 
we apply the Lagrange-d'Alembert principle stating that the action functional is critical with respect to variations subject 
to the constraint. In our case it reads
\begin{equation}\label{LdA}
\delta \int_0^TL_{ \Theta _0,d}(q, \dot q) dt=0, \quad \text{for variations}\quad  \delta q \in \mathcal{S}  _{ \bar\sigma }(q)
\end{equation} 
vanishing at $t=0,T$, and with $\dot q \in \mathcal{S}  _{ \bar\sigma }(q)$.

In a similar way with the continuous case, the relation between the Lagrangian variables 
$q(t) \in \mathsf{D} _{\bar \sigma } ( \mathbb{M}  )$, $ \Theta _0\in \mathbb{R}  ^N $,  
and the Eulerian variables $ A(t) \in \mathfrak{d}_{\bar \sigma }( \mathbb{M}  )$, $ \Theta (t)\in \mathbb{R} ^N $, is given by the formulas
\begin{equation}\label{formulas_Lagr_Euler} 
A(t)=\dot  q (t)q(t)^{-1} \quad\text{and}\quad  \Theta (t)= q(t)\Theta _0.
\end{equation} 
The discrete Lagrangian $L_{ \Theta _0,d}$ is assumed to have the same right-invariance as its continuous counterpart in \eqref{HP}, hence it can be exclusively written in terms of the Eulerian variables in \eqref{formulas_Lagr_Euler}. We thus get the reduced Lagrangian
\[
\ell_d=\ell_d(A, \Theta ):\mathfrak{d}_{\bar \sigma  } (\mathbb{M}) \times \mathbb{R}  ^N \rightarrow \mathbb{R}.
\]

The Eulerian version of the Lagrange-d'Alembert principle \eqref{LdA} is found to be
\begin{equation}\label{nonholonomic_variational_principle}
\delta \int_0^T\ell_d(A, \Theta)  dt=0, \quad \text{for variations} \quad   \delta A=\partial_tY+[Y,A],\quad \delta \Theta =Y \Theta ,
\end{equation}
in which $A \in \mathcal{S} _{ \bar \sigma }$ and where $Y $ is an arbitrary time 
dependent matrix in $\mathcal{S}_{ \bar \sigma }$ vanishing at the endpoints. The expressions for 
the variations $ \delta A$ and $ \delta \Theta $ in \eqref{nonholonomic_variational_principle} are 
obtained by using the two relations \eqref{formulas_Lagr_Euler}. In particular, we have $Y= \delta q q ^{-1} $, 
which is therefore an arbitrary time dependent matrix in $\mathcal{S} _{ \bar \sigma }$ vanishing at $t=0,T$. 
The principle \eqref{nonholonomic_variational_principle} is a nonholonomic version of the Euler-Poincar\'e variational principle.

\medskip

In order to derive the equations associated to the principle \eqref{nonholonomic_variational_principle}, we need to 
introduce appropriate dual spaces to the Lie algebra $\mathfrak{d}_{\bar \sigma  } (\mathbb{M})$ and the space of 
discrete functions $ \Omega _d ^0 ( \mathbb{M}  )= \mathbb{R}  ^N $.
To this end, we recall from \cite{GaMuPaMaDe2010} that in the context of discrete diffeomorphism groups, a discrete one-form on $\mathbb{M}$ is identified with a skew-symmetric $N \times N$ matrix. 
The space of discrete one-forms is denoted by $\Omega^1_d(\mathbb{M})$. The discrete exterior derivative of a discrete function 
$F \in  \Omega _d ^0 ( \mathbb{M}  )$ is the discrete one-form $\mbox{d}F$ given by
\[
(\mbox{d}F)_{ij}:=F_i-F_j.
\]
Then, given a strictly positive function $ \bar \sigma ( \mathbf{x} )>0$, the discrete version of the $L ^2 $ pairing \eqref{pairing_vector} is defined as
\begin{equation}\label{discrete_L2_pairing}
\left\langle \!\left\langle K, A \right\rangle \!\right\rangle _{ \bar \sigma }:= \operatorname{Tr} \left( K^\mathsf{T}\Omega ^{\bar \sigma }A \right) ,\quad K\in\Omega^1_d(\mathbb{M}),\quad A\in\mathfrak{d}_{\bar\sigma }(\mathbb{M}).
\end{equation} 
By repeating the arguments of Theorem 2.4 of \cite{GaMuPaMaDe2010} for the discrete $L^2$ pairing \eqref{discrete_L2_pairing} with weight $\bar \sigma $, we get the identification
\begin{equation}\label{identification}
\mathfrak{d}_{ \bar \sigma  }(\mathbb{M})^*\simeq \Omega^1_d(\mathbb{M})/\mbox{d}\Omega_d^0(\mathbb{M}),
\end{equation} 
which is the discrete analogue of the identification in Lemma \ref{Lemma2}.

Concerning functions, the discrete analogue of the pairing \eqref{pairing_function} is given by
\begin{equation}\label{discrete_L2_pairing_functions}
\left\langle F, G\right\rangle _{ \bar \sigma }:=  F^\mathsf{T}\Omega ^{\bar \sigma }G= \sum_{i=1}^NF_i\Omega_i ^{\bar \sigma }G_i,\quad F,G \in \mathbb{R}  ^N.
\end{equation} 

A direct application of the principle \eqref{nonholonomic_variational_principle} yields the following result.

\begin{proposition} A curve $(A(t), \Theta (t)) \in \mathfrak{d} _ {\bar \sigma }( \mathbb{M}  ) \times \mathbb{R}  ^N $ is critical for the 
principle \eqref{nonholonomic_variational_principle} if and only if there exists a discrete function $P \in \mathbb{R}  ^N $ such that the 
following equation holds
\begin{equation}\label{Discrete_EP} 
 \frac{d}{dt} \frac{\delta \ell_d}{\delta A _{ij} }
+ \left( \left[\frac{\delta \ell_d}{\delta A} \Omega ^{\bar \sigma }, A \right] (\Omega ^{\bar \sigma })^{-1} \right)  _{ij} 
- \frac{1}{2}  \left( \frac{\delta \ell_d}{\delta \Theta  _i}+\frac{\delta \ell_d}{\delta \Theta  _j} \right)  ( \Theta _j - \Theta _i ) + (P_i-P_j)=0,
\end{equation} 
for all $ j \in N(i)$, where the functional derivatives $\frac{\delta \ell_d}{\delta A} $ and $\frac{\delta \ell_d}{\delta \Theta } $ 
are computed with respect to the pairings \eqref{discrete_L2_pairing} and \eqref{discrete_L2_pairing_functions}. 
\end{proposition}

Equation \eqref{Discrete_EP} yields a structure-preserving spatial discretization of the Euler-Poincar\'e equation \eqref{EP_general} on the mesh $ \mathbb{M}$.
For the anelastic, resp., pseudo-incompressible equations, we will choose $\bar \sigma =\bar \rho $, resp., $\bar \sigma = \bar \rho \bar \theta $ and use in \eqref{Discrete_EP} suitable approximations
\begin{equation}\label{Lagrangians} 
\ell_d=\ell_d(A, \Theta ):\mathfrak{d}_{\bar \rho } (\mathbb{M}) \times \mathbb{R}  ^N \rightarrow \mathbb{R},  \quad \text{resp.} \quad  \ell_d=\ell_d(A, \Theta ): \mathfrak{d}_{\bar \rho \bar \theta } (\mathbb{M}) \times \mathbb{R}  ^N  \rightarrow \mathbb{R},
\end{equation} 
of the Lagrangians \eqref{Lagrangian_anelastic}, resp., \eqref{Lagrangian_pseudo_incomp}.

\paragraph{Anelastic system.} The discrete Lagrangian associated to \eqref{Lagrangian_anelastic} is 
\begin{equation}\label{l_d_anel} 
\ell_d(A, \Theta )= \frac{1}{2} \left\langle \!\left\langle A ^\flat , A\right\rangle \!\right\rangle _{\bar \rho } - c_p\left\langle \bar \Pi , \Theta \right\rangle _{\bar \rho }
    = \frac{1}{2} \sum_{ij} A _{ij} ^\flat A _{ij} \Omega ^{\bar \rho } _i - c _p \sum_i \bar \Pi _i \Theta_i\Omega ^{\bar \rho } _i ,
\end{equation} 
where the first, resp., the second duality pairing is given in \eqref{discrete_L2_pairing}, resp., \eqref{discrete_L2_pairing_functions}, 
and $\bar \Pi \in \mathbb{R}  ^N $ is a discretization of the reference value $\bar \pi (z)$ of the Exner pressure.
The first term in \eqref{l_d_anel} is the discretization of the kinetic energy associated to a given 
Riemannian metric on $ \mathcal{D} $ and is based on a suitable flat operator $A \in \mathcal{S} _{\bar \rho } \mapsto A ^\flat \in \Omega _d ^1 ( \mathbb{M}  )$ 
associated to the mesh $ \mathbb{M}$, see \cite{PaMuToKaMaDe2010}.

The functional derivatives of $\ell_d$ with respect to the pairings $ \left\langle \! \left\langle \,, \right\rangle \! \right\rangle _{ \bar \rho }$ 
and $ \left\langle \,, \right\rangle _{ \bar \rho }$ are, respectively, 
\[
\frac{\delta \ell_d}{\delta A_{ij} }= A_{ij} ^\flat \quad\text{and} \quad\frac{\delta \ell}{\delta \Theta _i }= - c _p \bar \Pi _i.
\]
Using them in \eqref{Discrete_EP} with $\bar \sigma =\bar \rho $, we get the structure-preserving spatial discretization of the anelastic system on the mesh $ \mathbb{M}  $ as
\begin{equation}\label{discrete_anel_general} 
\frac{d}{dt} A _{ij}^\flat +[A ^\flat \Omega ^{\bar \rho }, A ]_{ij} \frac{1}{\Omega ^{\bar \rho } _j} + c_p \frac{ \bar \Pi _i + \bar \Pi _j}{2} ( \Theta _j -\Theta _i) 
=- (P_i-P_j), \quad \text{for all $ j \in N(i)$}.
\end{equation} 

\paragraph{Pseudo-incompressible system.} The discrete Lagrangian associated to \eqref{Lagrangian_pseudo_incomp} is\begin{equation}\label{l_d_pi}
\ell_d(A, \Theta )= \frac{1}{2} \sum_{ij} \frac{1}{ \Theta _i }A _{ij} ^\flat A _{ij} \Omega ^{\bar \rho \bar \theta  } _i 
- g\sum_i \frac{1}{ \Theta _i }Z_i   \Omega ^{\bar \rho \bar \theta  } _i, 
\end{equation} 
in which $Z \in \mathbb{R}  ^N $ is a discretization of the height $z$. Note that we are now using the volumes 
$ \Omega ^{\bar \rho \bar \theta  }_i $.

The functional derivatives of $\ell_d$ with respect to the pairings $ \left\langle \! \left\langle \,, \right\rangle \! \right\rangle_{ \bar \rho \bar \theta }$ 
and  $ \left\langle \,, \right\rangle _{ \bar \rho \bar \theta }$  are, respectively,
\[
\frac{\delta \ell_d}{\delta A_{ij} }= \frac{1}{2} \left( \frac{1}{ \Theta _i }+\frac{1}{ \Theta _j } \right) A_{ij} ^\flat=: M_{ij}\ \  
\text{and}\ \ \frac{\delta \ell}{\delta \Theta _i }=  \frac{1}{ \Theta_i ^2 }\Big( gZ_i -k _i  \Big) , \ \  k_i := \frac{1}{2} \sum_j A _{ij} ^\flat A _{ij}.
\]
Using them in \eqref{Discrete_EP} with $\bar \sigma =\bar \rho \bar \theta $, we get the structure-preserving spatial discretization of the 
pseudo-incompressible system on the mesh $ \mathbb{M}  $ as
\begin{equation}\label{discrete_pi_general} 
\begin{split}
\frac{d}{dt} M_{ij} +[M \Omega ^{\bar \rho  \bar \theta }, A ]_{ij} \frac{1}{\Omega ^{\bar \rho  \bar \theta } _j} 
- \frac{1}{2}  \Big( \frac{gZ _i- k_i  }{ \Theta _i ^2 }  & +  \frac{gZ _j - k _j }{ \Theta _j ^2 } \Big) (\Theta_j -\Theta _i ) \\ & =- (P_i-P_j),  \quad \text{for all $ j \in N(i)$}.
\end{split}
\end{equation}

\section{Variational integrator on irregular simplicial meshes}
\label{sec_numerical_schemes}

\begin{figure}[t]
\centering
{\includegraphics[width=4.3in]{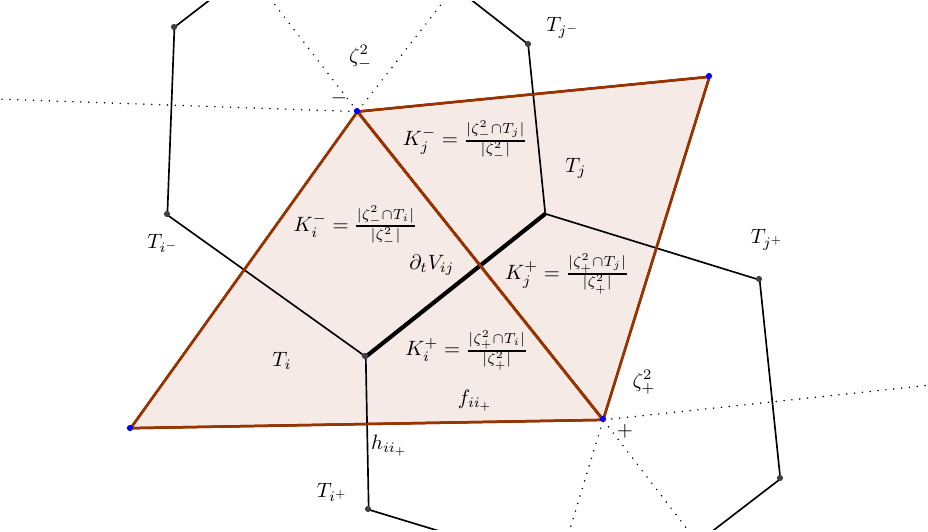}}
\caption{Notation and indexing conventions for the 2D simplicial mesh.}
\label{fig_Notation}
\end{figure}

In this section, we shall use the general results of \S\ref{VD}, valid for any kind of reasonable (i.e. non-degenerated) meshes, to deduce the variational discretization on 2D simplicial meshes.
On such meshes, we adopt the following notations (cf. Figure~\ref{fig_Notation}):
\begin{align*} 
f_{ij}:&=  \;\;\text{length of the primal edge, located between triangle $i$ and triangle $j$}; \\
h_{ij}:&=  \;\;\text{length of the dual edge that connect the circumcenters of triangle $i$ and triangle $j$};\\
\Omega_i:&=\;\; \text{area of the primal simplex (triangle) $T_i$}.
\end{align*} 
The flat operator on a 2D simplicial mesh is defined by the following two conditions, see \cite{PaMuToKaMaDe2010},
\begin{equation}\label{def_flat}
\begin{aligned} 
&A ^\flat _{ij}= 2\Omega_i\frac{h_{ij} }{f _{ij} }A _{ij}, \quad \text{for $j \in N(i)$},\\
&A ^\flat _{ij}+ A ^\flat _{jk}+ A^\flat _{ki}= K_j ^{e} \left\langle\omega (A ^\flat ),\zeta_{e} ^2 \right\rangle, \quad \text{for $i,k \in N(j)$, $k \notin N(i)$},
\end{aligned}
\end{equation}  
in which $e$ denotes the node common to triangles $T_i, T_j, T_k$ and $ \zeta ^2 _e $ denotes the dual cell to $e$. In \eqref{def_flat}, the vorticity $ \omega (K)$ of a discrete one-form $K \in \Omega ^1 _d ( \mathbb{M}  )$ is defined by
\[
\left\langle\omega (K),\zeta^2 _e\right\rangle:=  \sum_{ \zeta ^1 _{mn}\in \partial \zeta ^2  _e } K _{mn},
\]
where the sum is taken over the dual edges in the boundary $ \partial \zeta _e ^2 $
counterclockwise around node $e$. The constant $K _j ^e $ is defined as
\[
K^e_j:=\frac{|\zeta^2_{e}\cap T_j|}{|\zeta^2_{e}|},
\]
where $| \zeta _e ^2 |$ and $|\zeta ^2 _e\cap T_j|$ denote, respectively, the areas of $ \zeta ^2  _e$ and $ \zeta _e \cap T_j$. 
Note that the matrix $ A ^\flat $ defined in \eqref{def_flat} is skew-symmetric, hence $ A ^\flat \in \Omega _d ^1 (\mathbb{M}  )$.

%\todo{Here we need a picture of the mesh with the positions of $i$, $i_-$, $j$, $j_-$, .....}

\paragraph{Boussinesq flow.}
Variational discretization of the Boussinesq fluid on regular Cartesian grids has been carried out in 
\cite{DeGaGBZe2014}. Here we shall derive from \eqref{Discrete_EP} the variational scheme on irregular 2D simplicial grids.
Recall that in this case $\operatorname{div}_ \mu ( \mathbf{u} )=0$ and that the Boussinesq Lagrangian is given by
\[
\ell( \mathbf{u} , b)=\int_ \mathcal{D} \Big( \frac{1}{2} | \mathbf{u} |^2 +bz \Big) \mu . 
\]
The discrete Lagrangian is therefore chosen as $\ell_d: \mathfrak{d}( \mathbb{M}  ) \times \mathbb{R}  ^N \rightarrow \mathbb{R}  $,
\[
\ell_d(A, B )= \frac{1}{2} \left\langle \! \left\langle A ^\flat  , A  \right\rangle \!\right\rangle+ \left\langle B,Z \right\rangle,
\]
where $B \in \mathbb{R}  ^N $ is the discrete buoyancy and $Z $ is the discrete height function, i.e., $Z _i $ 
denotes the $z$-coordinate of the circumcenter of cell $C _i $.

Using the Boussinesq Lagrangian and the flat operator \eqref{def_flat}, the discrete Euler-Poincar\'e equation \eqref{Discrete_EP} yields
\begin{equation}\label{discrete_Boussinesq}
\left\{ 
\begin{array}{l}
\vspace{0.2cm}\partial _tA ^\flat _{ij} + \left\langle\omega(A ^\flat  ),\zeta^2_-\right\rangle\left(K_i^-A_{ii_-}+K^-_jA_{jj_-}\right)
                                        -\left\langle\omega (A ^\flat   ),\zeta ^2_+\right\rangle\left(K_i^+A_{ii_+}+K^+_jA_{jj_+}\right)\\
\vspace{0.2cm}\displaystyle\qquad\qquad \qquad \qquad\qquad \qquad  =\frac{Z _i + Z _j }{2} (B_j-B_i) + ( \tilde P_j-\tilde P_i), \quad \text{for all $ j\in N(i)$},\\
\displaystyle \partial _t B _i -\sum_{j\in N(i)}A _{ij} B _j=0,
\end{array} \right.
\end{equation}  
where $\Omega _i A_{ij}=- \Omega _j A_{ji}$, for all $i,j$, and $\sum_{j \in N(i)}A _{ij} =0$, for all $i$, and where $\tilde P$ is related to $P$ in \eqref{Discrete_EP} via
\begin{equation}\label{Ptilde} 
\tilde P _i = P _i + \sum_{k \in N(i)} A^\flat _{ik} A_{ik}.
\end{equation}
We note that the momentum equation \eqref{discrete_Boussinesq} corresponds to the discretization of the following form of the Boussinesq equation:
\begin{equation}\label{Boussinesq_vorticity} 
\partial _t \mathbf{u}^\flat + \mathbf{i} _\mathbf{u} \mathbf{d}\mathbf{u}^\flat 
 =- z \mathbf{d} b - \mathbf{d} \tilde p,
\end{equation} 
where, similarly to \eqref{Ptilde},  $\tilde p= \mathbf{i} _ \mathbf{u}\mathbf{u}^\flat +p$, with $p$ the pressure function arising in the Euler-Poincar\'e formulation \eqref{EP_general}. The form \eqref{Boussinesq_vorticity} is easily seen to be equivalent to the standard form \eqref{Boussinesq_standard_form} with $P_{\rm b} = zb+p + \frac{1}{2} | \mathbf{u} | ^2=zb+\tilde p - \frac{1}{2} | \mathbf{u} | ^2 $.

\paragraph{Anelastic flow.}
The continuous and discrete anelastic Lagrangians are given in \eqref{Lagrangian_anelastic} and \eqref{l_d_anel}. Recall that in this case $\operatorname{div}_ \mu  ( \bar \rho \mathbf{u}) =0$. 
The flat operator  \eqref{def_flat} has to be slightly modified in order to produce a skew-symmetric matrix, namely, 
we modify the first line in \eqref{def_flat} to
\begin{equation}\label{def_flat_new} 
A ^\flat := M^{(A)},\;\; \text{for the matrix $M$ defined by} \;\; M_{ij}:=2\Omega_i\frac{h_{ij} }{f _{ij} }A _{ij},
\end{equation} 
in which $( \cdot ) ^{(A)}$ denotes the skew-symmetric part. For the Boussinesq model,
this definition recovers \eqref{def_flat}, since the matrix $M$ is 
in this case already skew-symmetric.
One checks that this definition still satisfies the properties of a flat operator in \cite{PaMuToKaMaDe2010}.

The general discrete anelastic equations \eqref{discrete_anel_general} yield
\begin{equation}\label{discrete_anelastic_2D}
\left\{ 
\begin{array}{l}
\vspace{0.2cm}\partial _tA ^\flat _{ij} 
+ \left\langle\omega(A ^\flat  ),\zeta^2_-\right\rangle\left(K_i^-A_{ii_-}+K^-_jA_{jj_-}\right)
 -\left\langle\omega (A ^\flat   ),\zeta ^2_+\right\rangle\left(K_i^+A_{ii_+}+K^+_jA_{jj_+}\right)\\
\vspace{0.2cm}\displaystyle\qquad\qquad \qquad \qquad\qquad =-  c_p \frac{ \bar \Pi _i + \bar \Pi _j}{2} ( \Theta _j -\Theta _i)+ ( \tilde P_j-\tilde P_i), 
\quad \text{for all $ j\in N(i)$}, \\
\displaystyle \partial _t \Theta  _i -\sum_{j\in N(i)}A _{ij} \Theta  _j=0,
\end{array} \right.
\end{equation}  
where  $\Omega _i ^{\bar \rho } A_{ij}=- \Omega _j ^{\bar \rho }A_{ji}$, for all $i,j$, and $\sum_{ j\in N(i)}A _{ij} =0$, for all $i$, and where $\tilde P$ is related to $P$ in \eqref{Discrete_EP} and \eqref{discrete_anel_general} as before via the formula \eqref{Ptilde}.
We note that the momentum equation \eqref{discrete_anelastic_2D} corresponds to the discretization of the following form of the anelastic equation:
\begin{equation}\label{anelastic_vorticity} 
\partial _t \mathbf{u}^\flat + \mathbf{i} _\mathbf{u} \mathbf{d}\mathbf{u}^\flat 
 =c_p \bar \pi \mathbf{d} \theta  - \mathbf{d} \tilde p,
\end{equation} 
where, similarly to \eqref{Ptilde},  $\tilde p= \mathbf{i} _ \mathbf{u}\mathbf{u}^\flat +p$, with $p$ the pressure function arising in the Euler-Poincar\'e formulation \eqref{EP_general}. The form \eqref{anelastic_vorticity} was shown in \S\ref{VF} to be equivalent to the standard form \eqref{anelastic0}.

\paragraph{Pseudo-incompressible flow.}
The continuous and discrete pseudo-incompressible Lagrangians are given in \eqref{Lagrangian_pseudo_incomp} and \eqref{l_d_pi}. Recall that in this case $\operatorname{div}_ \mu  ( \bar \rho \bar \theta \mathbf{u}) =0$. 
We take the flat operator \eqref{def_flat} with the first line modified as in \eqref{def_flat_new}

The general discrete pseudo-incompressible equations \eqref{discrete_pi_general} yield
\begin{small}
\begin{equation}\label{discrete_rot_Euler_PI}
\left\{ 
\begin{array}{l}
\vspace{0.2cm}\partial _tM_{ij} 
+ \left\langle\omega( M ),\zeta^2_-\right\rangle\left(K_i^-A_{ii_-}+K^-_jA_{jj_-}\right)
-\left\langle\omega ( M),\zeta ^2_+\right\rangle\left(K_i^+A_{ii_+}+K^+_jA_{jj_+}\right)\\
\vspace{0.2cm}\displaystyle\qquad\qquad\qquad 
=\frac{1}{2} \left( \frac{gZ _i- k_i  }{ \Theta _i ^2 } +\frac{gZ _j - k _j }{ \Theta _j ^2 } \right) (\Theta_j -\Theta _i ) + ( \tilde P_j-\tilde P_i), \quad \text{for all $ j\in N(i)$}\,, \\
\vspace{0.2cm}\displaystyle M _{ij} =\frac{1}{2} \left( \frac{1}{ \Theta _i }+\frac{1}{ \Theta _j } \right) A_{ij} ^\flat,\\
\displaystyle \partial _t \Theta  _i -\sum_{j\in N(i)}A _{ij} \Theta  _j=0,
\end{array} \right.
\end{equation}
\end{small}  
\!\!where  $\Omega _i ^{\bar \rho \bar \theta } A_{ij}=- \Omega _j ^{\bar \rho \bar \theta }A_{ji}$, for all $i,j$, and $\sum_{ \in N(i)}A _{ij} =0$, for all $i$, and where $\tilde P$ is related to $P$ in \eqref{Discrete_EP} and \eqref{discrete_pi_general} by the fomula
\begin{equation}\label{Ptilde_pi} 
\tilde P _i = P _i + \sum_{k \in N(i)} M _{ik} A_{ik} .
\end{equation} 
We note that the momentum equation \eqref{discrete_rot_Euler_PI} corresponds to the discretization of the 
following form of the pseudo-incompressible equation:
\begin{equation}\label{pi_vorticity} 
\partial _t \left( \frac{1}{ \theta } \mathbf{u}^\flat \right) + \frac{1}{ \theta } \mathbf{i} _\mathbf{u} \mathbf{d}\mathbf{u}^\flat 
 =-\frac{1}{\theta ^2 }(gz- \frac{1}{2} | \mathbf{u} | ^2 )\mathbf{d} \theta    - \mathbf{d} \tilde p,
\end{equation} 
where, similarly to \eqref{Ptilde_pi},  $\tilde p=\frac{1}{\theta } \mathbf{i} _ \mathbf{u}\mathbf{u}^\flat +p$, with $p$ the 
pressure function arising in the Euler-Poincar\'e formulation \eqref{EP_general}. The form \eqref{pi_vorticity} was shown in \S\ref{VF} 
to be equivalent to the standard form \eqref{pseudo_incompressible}.

We present in Table \ref{table} a parallel between the continuous and discrete variational formulations for the three models.

\paragraph{Time integration.} Since the spatial discretization has been realized in a structure-preserving way, a corresponding temporal variational discretization follows by applying the general discrete (in time) Euler-Poincar\'e approach,
as it has be done in \cite{GaMuPaMaDe2010}, \cite{DeGaGBZe2014} to which we refer for a detailed treatment. This approach is based on the use of the Cayley transform, a local approximant of the exponential map.
For the general discrete Euler-Poincar\'e system \eqref{Discrete_EP} and for a given time step $\Delta t$, it results in the following scheme
\begin{align*}
&\frac{1}{\Delta t} \bigg( \frac{\delta \ell_d}{\delta A ^k  _{ij} }-  \frac{\delta \ell_d}{\delta A ^{k-1}  _{ij} }\bigg) 
+ \frac{1}{2} \left( \left[\frac{\delta \ell_d}{\delta A^k } \Omega ^{\bar \sigma }, A ^k  \right] (\Omega ^{\bar \sigma })^{-1}+\left[\frac{\delta \ell_d}{\delta A^{k-1}} \Omega ^{\bar \sigma }, A^{k-1} \right] (\Omega ^{\bar \sigma })^{-1} \right)  _{ij} \\
& \qquad \qquad \qquad \qquad \qquad \qquad\qquad\qquad 
- \frac{1}{2}  \bigg( \frac{\delta \ell_d}{\delta \Theta  ^k _i}+\frac{\delta \ell_d}{\delta \Theta  ^k _j} \bigg)  ( \Theta ^k _j - \Theta ^k _i ) + (P ^k _i-P ^k _j)=0,
\end{align*} 
where $(A _{ij} ^{k-1}, \Theta ^{k-1}_i)$ and $(A _{ij}^k , \Theta_i ^k )$ are the values at the consecutive time steps $k-1$ and $k$.

\begin{small} 
\begin{table}[h!]\centering 
\centering
\begin{tabular}{|c | c |} 
\hline
{\bf Continuous diffeomorphisms} & {\bf Discrete diffeomorphisms}\\ [0.45ex] 
\hline
Boussinesq: \;\;$\operatorname{Diff}_\mu (M)$ & Boussinesq: \;\;$\mathsf{D}( \mathbb{M}  )$ \\
Anelastic: \;\;$\operatorname{Diff} _{\bar{ \rho  }\mu}(M)$ & Anelastic: \;\;$\mathsf{D}_{\bar \rho }( \mathbb{M}  )$ \\
Pseudo-incompressible: \;\;$\operatorname{Diff} _{\bar{ \rho  }\bar\theta \mu} (M)$ &Pseudo-incompressible: \;\;$\mathsf{D}_{\bar \rho \bar\theta }( \mathbb{M}  )$ \\[0.45ex] 
\hline\hline
{\bf Lie algebras} & {\bf Discrete Lie algebras}\\ [0.45ex] 
\hline
$ \mathfrak{X}_ \mu (M),\;\; \mathfrak{X}_ {\bar{ \rho  }\mu} (M), \;\;\mathfrak{X}_{ \bar{ \rho  }\bar \theta \mu }(M) $ & $ \mathfrak{d} ( \mathbb{M}  ), \;\;\mathfrak{d} _ {\bar{ \rho  } } ( \mathbb{M}  ), \;\;\mathfrak{d} _{ \bar{ \rho  }\bar \theta   }( \mathbb{M}  )$ \\[0.45ex] 
\hline\hline
{\bf Euler-Poincar\'e form} & {\bf Discrete Euler-Poincar\'e form}\\ [0.45ex] 
\hline
$\displaystyle\partial _t \frac{\delta \ell}{\delta \mathbf{u} }+\pounds _\mathbf{u} \frac{\delta \ell}{\delta \mathbf{u} } + \frac{\delta \ell}{\delta \theta } \mathbf{d} \theta = - \mathbf{d} p$, 
&Equation \eqref{Discrete_EP} \\ [0.95ex]
Common form for the three models &Common discrete form for the three models\\
 &Form independent of the mesh\\
[0.25ex]
\hline\hline
{\bf Expression corresponding to the} & {\bf Discrete form on 2D simplicial grids}\\ 
{\bf discrete form on 2D simplicial grids} & \\[0.45ex] 
\hline
Boussinesq: & Discrete Boussinesq:\\
$\displaystyle\partial _t \mathbf{u}^\flat + \mathbf{i} _\mathbf{u} \mathbf{d}\mathbf{u}^\flat 
 =- z \mathbf{d} b - \mathbf{d} \tilde p$ & Equation \eqref{discrete_Boussinesq}  \\[0.45ex]
Anelastic: & Discrete Anelastic: \\ [0.45ex]
$\displaystyle\partial _t \mathbf{u}^\flat + \mathbf{i} _\mathbf{u} \mathbf{d}\mathbf{u}^\flat 
 =c_p \bar \pi \mathbf{d} \theta  - \mathbf{d} \tilde p$ & Equation \eqref{discrete_anelastic_2D} \\[0.45ex]
Pseudo-incompressible: & Discrete Pseudo-incompressible: \\ [0.45ex]
$\displaystyle\partial _t \Big( \frac{1}{ \theta } \mathbf{u}^\flat \Big) + \frac{1}{ \theta } \mathbf{i} _\mathbf{u} \mathbf{d}\mathbf{u}^\flat 
 =-\frac{1}{\theta ^2 }\Big(gz- \frac{1}{2} | \mathbf{u} | ^2 \Big)\mathbf{d} \theta    - \mathbf{d} \tilde p$ & Equation \eqref{discrete_rot_Euler_PI} \\[0.25ex]
\hline
\end{tabular}
\vspace{1em}
\caption{Parallel between the continuous and discrete forms for the three models. 
Note that in the Euler-Poincar\'e form given in the sixth row of the first column, one has to compute 
the variational derivatives with respect to the three different weighted pairings in order to get the three models. 
The last row of the first column presents the continuous equations in a form that corresponds to the discrete 
forms obtained by variational discretization on 2D simplicial meshes. Note that these expressions are not the 
standard form of the models given in \eqref{Boussinesq_standard_form}, \eqref{anelastic}, \eqref{pseudo_incompressible}.} \label{table} 
\end{table}
\end{small}

  \color{black}

\section{Numerical tests}
\label{sec_numerical_tests}
 
  In this section we present preliminary numerical tests for the variational schemes. 
  We will focus on hydrostatic adjustment processes and make for each model a quantitative 
  evaluation of the discrete dispersion relation of the emitted internal gravity waves. %, for each of the three models.  
  The simulations are performed on a regular and an irregular triangular mesh.

  \paragraph{Description of the meshes.}
  The regular mesh consists of equilateral triangles of constant edge length $f = |f_{ij}|$,  
  where $f_{ij}, j = 1,2,3,$ denote the edges of triangle $T_i$ (cf. Section \ref{sec_numerical_schemes}).
  The distance between neighboring vertices in $x$-direction is given by $f_x:= f$ while the height of
  the triangles in $z$-direction is given by $f_z := \frac{\sqrt{3}}{{2}} f$. 
  Given a domain size of $L_x \times L_z$, in which $L_x$ and $L_z$ denote the domain's length in $x$- and 
  $z$-directions, respectively, the {\em mesh resolution}, denoted by $ 2\cdot N_x \times N_z$ 
  for $N_x :=  {L_x}/{f_x}$ and $N_y :=  {L_z }/{f_z}$, corresponds to the number of triangular cells.

  To construct the irregular mesh, we start from the regular one and randomly move the regularly distributed {\em internal} vertices 
  -- i.e. vertices that do {\em not} belong to boundary cells %(triangles that describe the outer boundaries of the computational domain) 
  -- from point ${\bf x_i} = (x_i,z_i)$ to  ${\bf x_i} + \delta {\bf x_i}$ 
  within the bounds $ |\delta x_i| < c\cdot f_x \cdot r$ and $ |\delta z_i| < c\cdot f_z \cdot r$, for a positive constant $c$
  and some random number $r \in [-0.5, 0.5 ]$. 
  Although not necessary, we leave the boundary triangles regular as this 
  eases the implementation.
  
  The distortion of the irregular mesh can be quantified using a grid quality measure 
  introduced in \cite{Baueretal2014} that measures the distortion of the dual cells:
  $\Delta h({\bf x})  := \frac{\max_j \ h_{ij} }{\min_j \ h_{ij}  } \ $,  
  in which $h_{ij}$ is the length of dual edge $j$ of dual cell $\zeta^2_i$ that contains point $\bf x$. 
  High values of $\Delta h$ indicate strongly deformed cells. 
  For our studies we use $c =  0.2$ which leads to a mesh with $\max_{{\bf x} \in \Omega}\Delta h({\bf x}) \approx 7$ 
  indicating strongly deformed dual mesh cells. 
%   This allows us to illustrate that our variational schemes 
%   are stable on such meshes while conserving mass and total energy to a high degree.

 \begin{figure}[t] \centering
  \begin{tabular}{c}    
  {\includegraphics[width=2.4in]{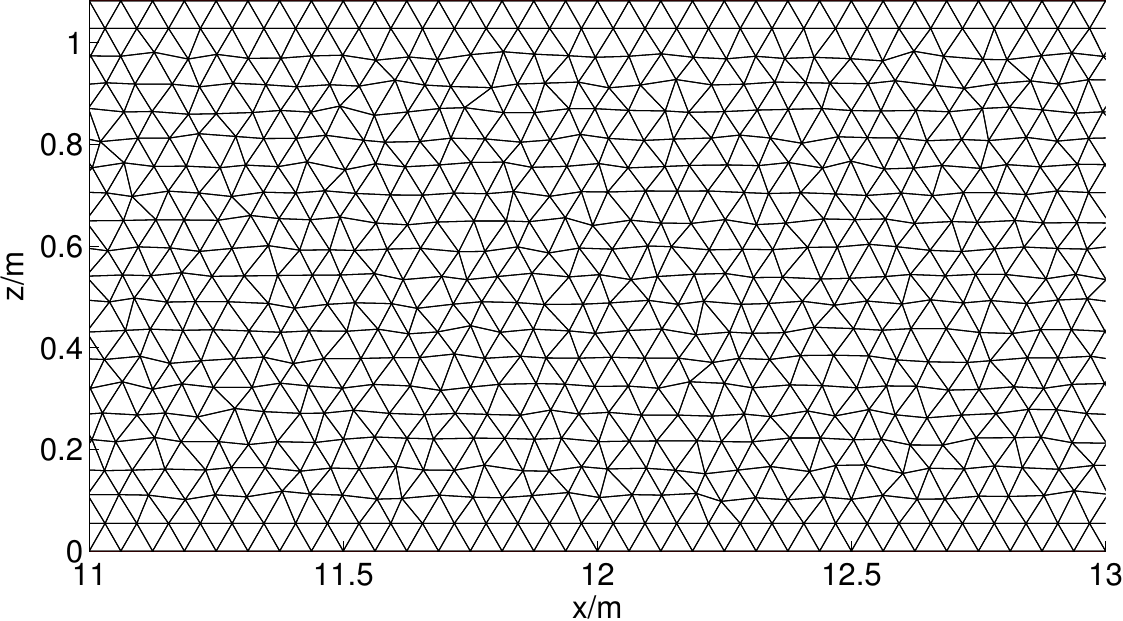}}
  \end{tabular}
    \caption{Section of central part of the irregular mesh with $\max_{{\bf x} \in \Omega}\Delta h({\bf x}) \approx 7$
    for a resolution of $2\cdot 384 \times 20$ triangular cells.}                                                                                              
    \label{fig:irr_grid}
  \end{figure}
 
 \color{black}
 
 We use a computation domain of dimension 
 $(x,z) \in \mathcal{D} =[0,L_x] \times [0, L_z], L_x = 24\,{\rm m}, L_z = 1\,{\rm m}$, while imposing %\incl{L_z = 1.0825}
 periodic boundary conditions in $x$-direction and free-slip boundary conditions at the 
 upper and lower boundaries of the domain. Both regular and irregular computational meshes have a 
 resolution of $2\cdot 384 \times 20$ triangular cells (cf. Figure~\ref{fig:irr_grid}).

 \paragraph{Description of the hydrostatic adjustment test case.}

  The derivations of Boussinesq, anelastic, and pseudo-incompressible models
  rely on the assumption of a vertically varying reference state that is in hydrostatic balance, 
  i.e. the gravitational and pressure terms compensate each other 
  (cf. Section \ref{sec_anelastic_pseudoincomp_systems}). 
  When out of equilibrium, the system tends to a balanced state by the so-called hydrostatic adjustment 
  process \cite{Lamb1932} by emitting internal gravity waves. 
  
  Applying this test 
  case, we study the schemes' dynamical behavior, long term energy and mass conservation properties, 
  and their discrete dispersion relations. 
  We initialize the Boussinesq scheme as in \cite{DeGaGBZe2014}, and adapt the therein suggested 
  test case to suit also for the anelastic and pseudo-incompressible schemes. This will allow us 
  to compare quantitatively the simulation results of our schemes with each other and with 
  those of \cite{DeGaGBZe2014}.
 
\begin{figure}[t] \centering
 \begin{tabular}{cc}    
 \hspace{-0.2cm}{\includegraphics[width=3in]{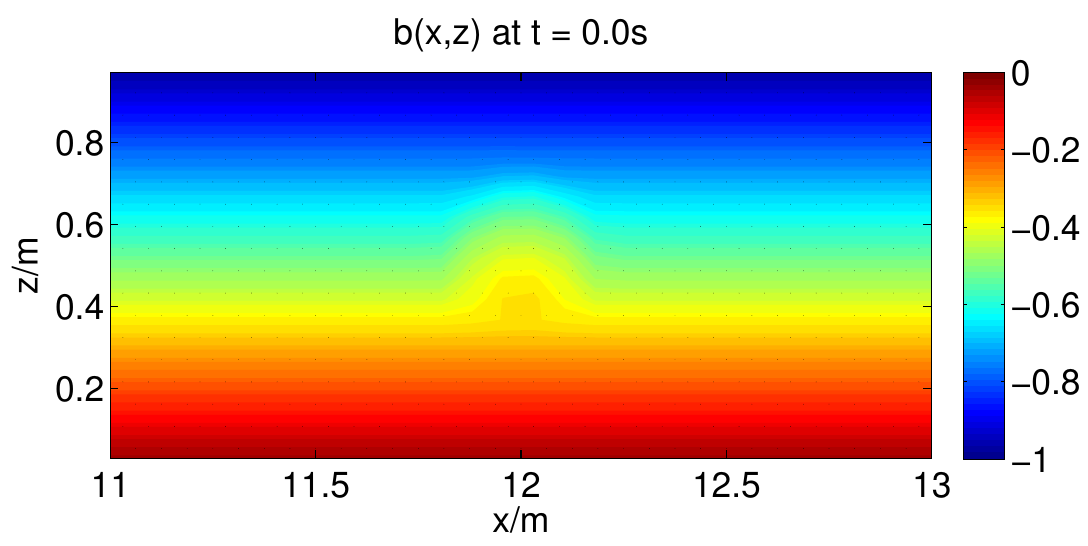}} & 
 \hspace{-0.0cm}{\includegraphics[width=3in]{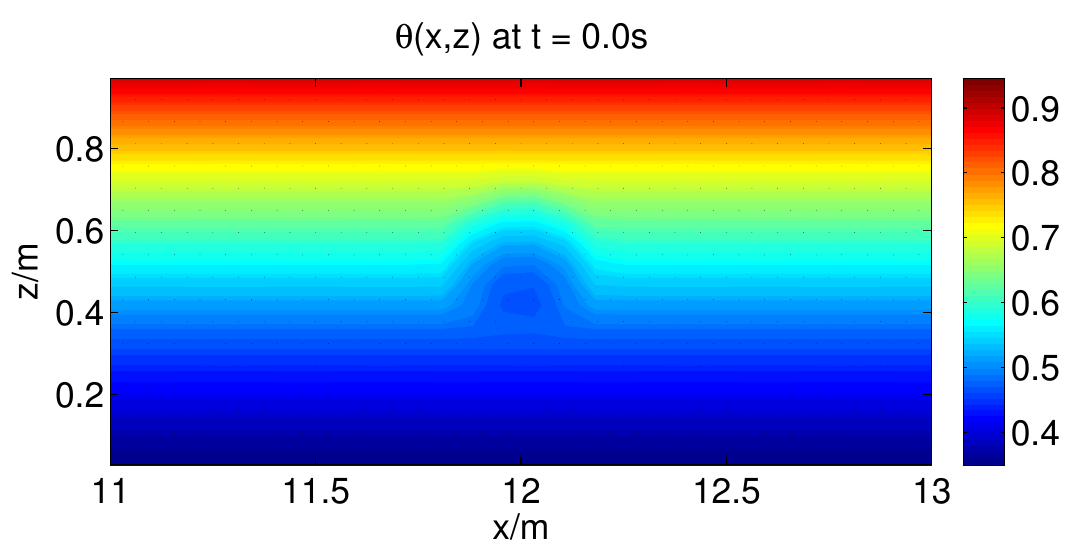}}
 \end{tabular}
  \caption{Initialization of the Boussinesq scheme by the buoyancy field $b(x,z,0)$, shown left. 
  Initialization of the anelastic and pseudo-incompressible schemes
  by the potential temperature field $\theta(x,z,0)$, shown right. 
  }                                                                                             
  \label{fig:initial_fqcy_dynamics}
 \end{figure}

  \paragraph{Initialization.} Analogously to \cite{DeGaGBZe2014}, we initialize the Boussinesq scheme 
  on the basis of a hydrostatic equilibrium,
  given by $u_{\rm eq}(x,z) = w_{\rm eq}(x,z) = 0$ and $b_{\rm eq}(x,z) = - N_b^2 z = : \bar b(z)$, on which 
  at $t=0$ a localized positive buoyancy disturbance $\tilde b(x,z)$ with compact support is superimposed.
  Hence, the initial buoyancy field $b(x,z,0) = \bar b(z) + \tilde b(x,z)$ with 
  Brunt-V\"as\"al\"a frequency $ N_b = 1/{\rm s}$ is given by the function 
  \begin{equation} \label{equ_init_oszill_Boussines}
   b(x,z,0) = N_b^2
         \begin{cases} 
           -z + \beta_b e^{\left(\frac{-r_0^2}{r_0^2 - r^2}\right)}  &\mbox{if } r < r_0  , \ \  r^2 = (x - \frac{L_x}{2})^2 + (z - \frac{L_z}{2})^2 ,\\
           -z                                          &\mbox{if } r \geq r_0  ,
         \end{cases} 
  \end{equation}
  with parameters $r_0 = 0.2 \cdot L_z$ and $\beta_b =  0.3 \cdot L_z$. %, $r^2 = (x - \frac{L_x}{2})^2 + (z - \frac{L_z}{2})^2$. 
  %As initial wind field we set $u_0(x,z) = v_0(x,z) = 0$. 
  Note that $[z] = \text{m}$, 
  hence the choice of $N_b = 1/{\rm s}$ suggests further to set $g = 1 \, \rm{ m /s^2}$ and $\theta_0 =1\, {\rm K}$. 
  Given these analytical functions, the discrete function $B = \{B_i|\ \text{for all triangles}\ T_i\}$ 
  is obtained by setting $B_i(0) = b(x_i,z_i,0)$ for all triangles $T_i$ with cell centers at position $(x_i,z_i)$
  (cf. Figure~\ref{fig:initial_fqcy_dynamics}).

  For the anelastic and pseudo-incompressible schemes, 
  we aim for an initialization that produces results comparable to the 
  Boussinesq scheme and that meets the requirements of constant $N$ and $\sigma_a$ or $\sigma_{\rm pi}$
  (discussed later in more detail). %in equations \eqref{N_sigma}. 
%   \werner{This will allow us to compare the 
%   results obtained by the anelastic and PI schemes with those of the Boussinesq scheme 
%   and with those in \cite{DeGaGBZe2014}.
%   }
  To this end, the hydrostatic equilibrium is set up by 
  $u_{\rm eq}(x,z) = w_{\rm eq}(x,z) = 0$ and a reference state $\bar \theta (z) = e^{z+c}$ with constant $c$,
  on which at $t=0$ a negative potential temperature perturbation $\tilde \theta(x,z)$ % \approx - \tilde b(x,z)$
  is superimposed. The initial potential temperature field 
  $\theta(x,z,0) = \bar \theta(z) + \tilde \theta(x,z)$ 
  is hence given by
  \begin{equation}\label{equ_init_oszill_AnPi}
  \theta(x,z,0) = 
         \begin{cases} 
           e^{z+c} - \beta_a e^{\left(\frac{-r_0^2}{r_0^2 - r^2}\right)}  &\mbox{if } r < r_0 , \quad r^2 = (x - \frac{L_x}{2})^2 + (z - \frac{L_z}{2})^2, \\
           e^{z+c}                                         &\mbox{if } r \geq r_0,
         \end{cases} 
 \end{equation}
 with parameters $r_0 = 0.2 \cdot L_z$ and $\beta_a =  0.2 \cdot L_z$. 
 To obtain a potential temperature field with comparable magnitude (in the order of $\theta_0 =1\,$K) 
 to the buoyancy field, we set $c = - L_z$ giving $0.4\,$K at the bottom and $1\,$K at the top 
 of the domain. 
 The choice of $\beta_a$ results in an oscillation 
 comparable in magnitude to the Boussinesq case. 
 For this $\theta$, the Brunt-V\"ais\"al\"a frequency is
 $N^2 = \frac{g}{\bar \theta }\frac{ d\bar\theta }{ dz} = 1/{\rm s}^2$, 
 where we set $g = 1\,{\rm m/s}^2$. 
 The requirement that $\sigma_a$ and $\sigma_{\rm pi}$ have to be constant 
 restricts our choice of the stratified density field $\bar \rho$ 
 to be either a constant or an exponential function; here we use 
 the profile $\bar \rho (z) = e^{- z}$ which mimics a realistic 
 stratification of the atmosphere.
 
 The initialization of the anelastic scheme requires, in addition, to define 
 a discrete Exner pressure $\bar \Pi$. % $\bar \pi(z)$. 
 The relation between $\bar \pi$ and $\bar \theta$, see \eqref{HB}, allows us to initialize 
 the Exner pressure by the potential temperature field via 
  \begin{equation}\label{equ_relationExnerPT_real}
   \bar \pi (z)  =   \frac{g}{c_p} e^ {- (z+c)} \, ,
  \end{equation}
 for any values of specific heat at constant pressure $c_p$ (here we set $c_p=1$),
 as it will cancel out in the anelastic equations.
 Given these functions, the discrete ones are obtained by setting 
 $\Theta_i(0) = \theta(x_i,z_i,0)$, $\bar \Pi_i = \bar \pi(x_i,z_i)$,  and $\bar \rho_i = \bar \rho (x_i,z_i)$ 
 for all triangles $T_i$ with cell centers at position $(x_i,z_i)$ (cf. Figure~\ref{fig:initial_fqcy_dynamics}).
  
 We integrate for a time interval of $100\,$s (in correspondence to \cite{DeGaGBZe2014}) 
 and use a fixed time step size of $\Delta t = 0.25\,$s for all schemes.

 \paragraph{Conserved quantities.}

 In Section~\ref{sec_anelastic_pseudoincomp_systems}, we considered soundproof models that provide 
 energy conserving approximations of the Euler equations. In the following we study if the 
 variational schemes conserve discrete versions of the associated total energies too.
 
 In the same vein, we study if discrete versions of mass are conserved quantities in time also.
 We note that mass conservation in the Boussinesq case is given by 
 \begin{equation}\label{equ_buyancymass_conservation}
  \frac{d}{dt} \int_ \mathcal{D}  b(x,z,t) d\mathbf{x} = 0.
 \end{equation}
 Being implicitly related to the density, we refer to this quantity, and the upcoming 
 similar ones for anelastic and pseudo-incompressible equations, generally as mass $M(t)$.
 For the anelastic equations, mass conservation is given by
  \begin{equation}\label{equ_anelasticymass_conservation}
  \frac{d}{dt} \int_{\mathcal{D}  }\bar \rho(z) \theta(x,z,t) d\mathbf{x} = 0 , 
  \end{equation}
  as $ \int_\mathcal{D} \bar \rho \partial_t \theta d\mathbf{x}  = - \int_\mathcal{D} \mathbf{d}  \theta \cdot \bar \rho \mathbf{u} d\mathbf{x} = 
  - \int_\mathcal{D} \operatorname{div}  (\bar \rho\mathbf{u}   \theta )  d\mathbf{x} 
  + \int_\mathcal{D}  \operatorname{div}  (\bar \rho \mathbf{u} ) \theta d\mathbf{x} =0$, 
  which follows by the anelastic constraint $ \operatorname{div}  (\bar \rho\mathbf{u} )=0$ and by the choice of boundary conditions, i.e. 
  ${\bf u} \cdot {\bf n} = 0$, on $\partial \mathcal{D} $.
  Following a similar argumentation, mass conservation for the pseudo-incompressible equations
  is given by
 \begin{equation}\label{equ_pseudomass_conservation}
  \frac{d}{dt} \int_ { \mathcal{D}}  \bar \rho(z) \bar \theta(z) \theta(x,z,t) d\mathbf{x}  = 0 . 
 \end{equation}

 \paragraph{Results on the dynamics.}

 \begin{figure}[t] \centering
 \begin{tabular}{cc}    
 \hspace*{-0.2cm}{\includegraphics[width=3in]{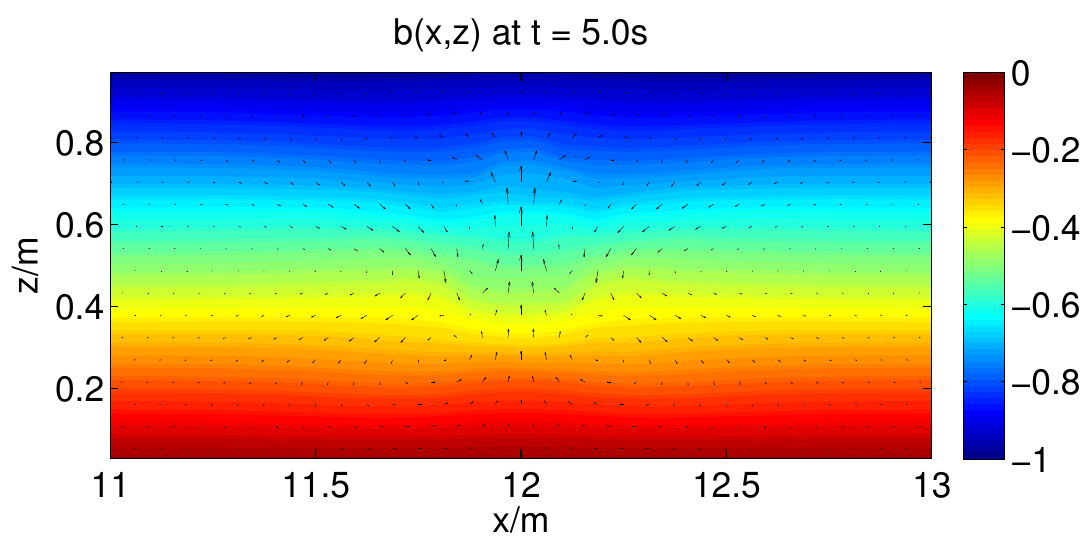}} & 
 \hspace*{-0.0cm}{\includegraphics[width=3in]{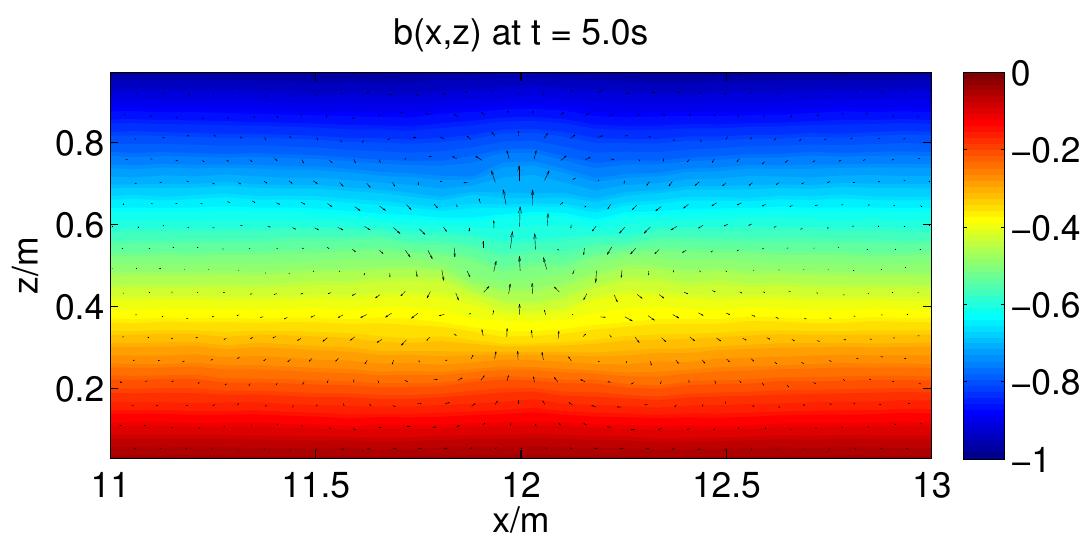}}\\ 
 \hspace*{-0.2cm}{\includegraphics[width=3in]{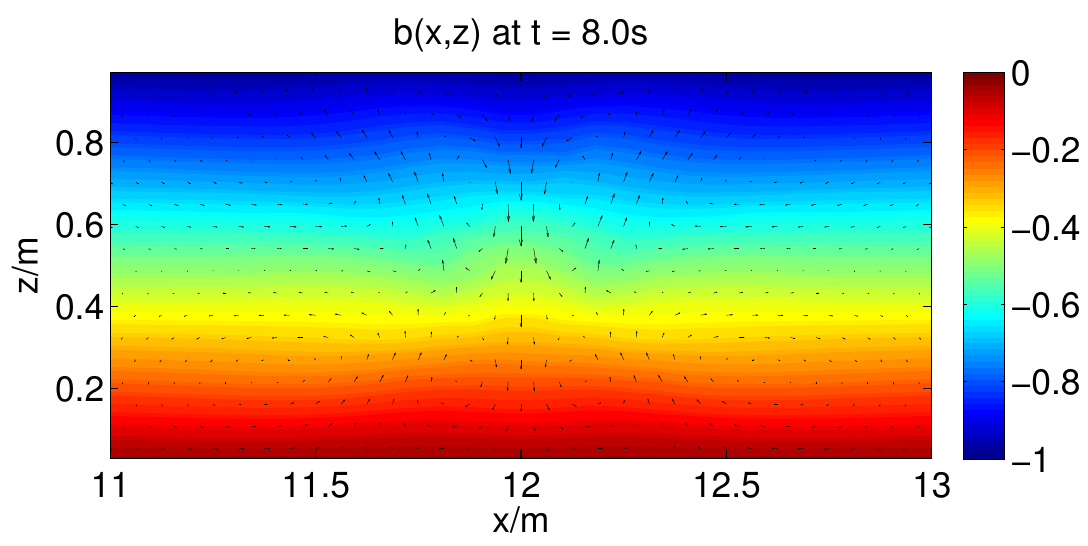}} & 
 \hspace*{-0.0cm}{\includegraphics[width=3in]{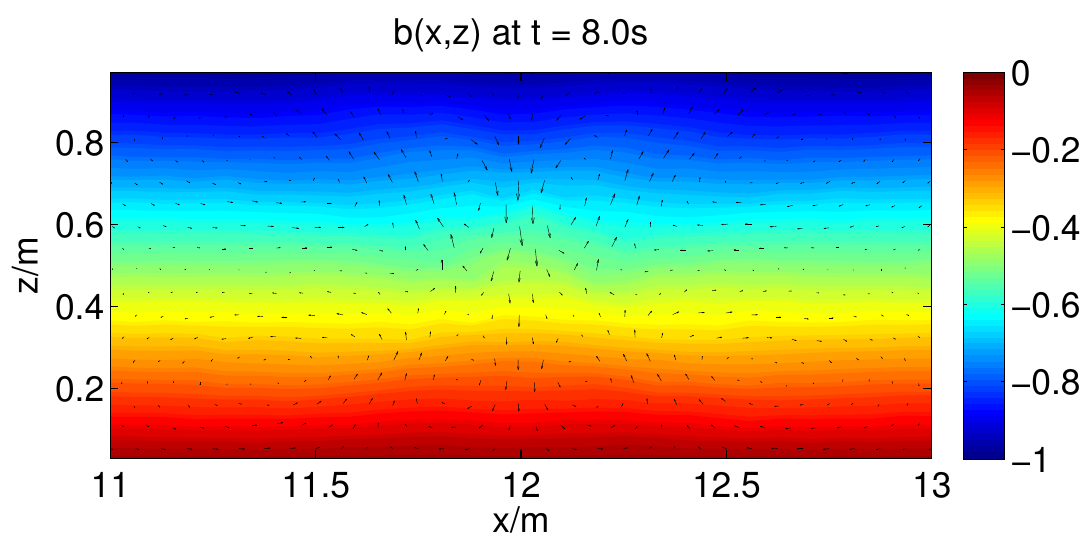}} 
 \end{tabular}
  \caption{Boussinesq scheme: snapshots of the wave propagation on the regular (left column) and 
  the irregular (right column) mesh.
  }                                                                                             
  \label{fig:byncy_dynamics}
 \end{figure}

 Before discussing the quantities of interest, let us first have a look at the general dynamical 
 behavior of the variational schemes. 
 Figure~\ref{fig:byncy_dynamics} shows snapshots at times $t = 5\,$s and $t = 8\,$s 
 of the buoyancy field $b(x,z,t)$ of the Boussinesq scheme
 for the central region $[11\,{\rm m},13\,{\rm m} ] \times [0,1\,{\rm m}]$ of the regular (left column) and 
 the irregular (right column) mesh. 
%   Figure~\ref{fig:byncy_dynamics} shows for regular (left column) and irregular (right column) meshes
%  snapshots of the buoyancy field $b(x,z,t)$ \werner{of the Boussinesq scheme} at times $t = 5\,$s and $t = 8\,$s for a region 
%  $[11\,{\rm m},13\,{\rm m} ] \times [0,1\,{\rm m}]$. 
 For these early times, before waves that are reflected by the boundaries reach the center,
 one clearly observes the internal gravity waves, caused by the buoyancy perturbation, 
 that propagate from the center along the channel in $x$-direction.
 Besides of small irregularities of the solutions on the irregular mesh, in particular visible at the velocity field 
 that is not completely symmetric with respect to the axis $x=10\,{\rm m}$, the results obtained
 using either the regular or the irregular mesh are very similar.

\begin{figure}[t] \centering
 \begin{tabular}{cc}
 \hspace*{-0.2cm}{\includegraphics[width=3in]{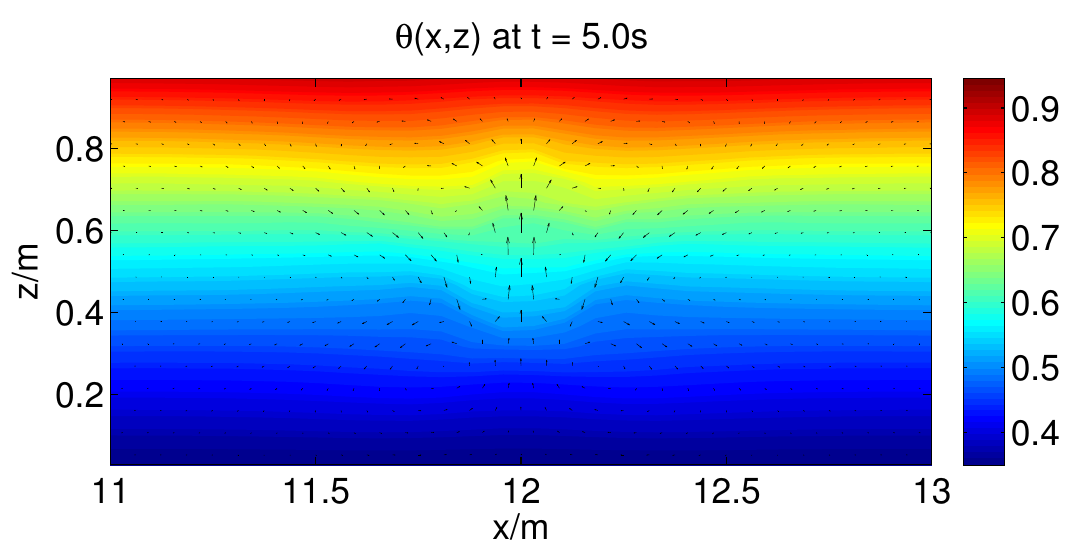}} & 
 \hspace*{-0.0cm}{\includegraphics[width=3in]{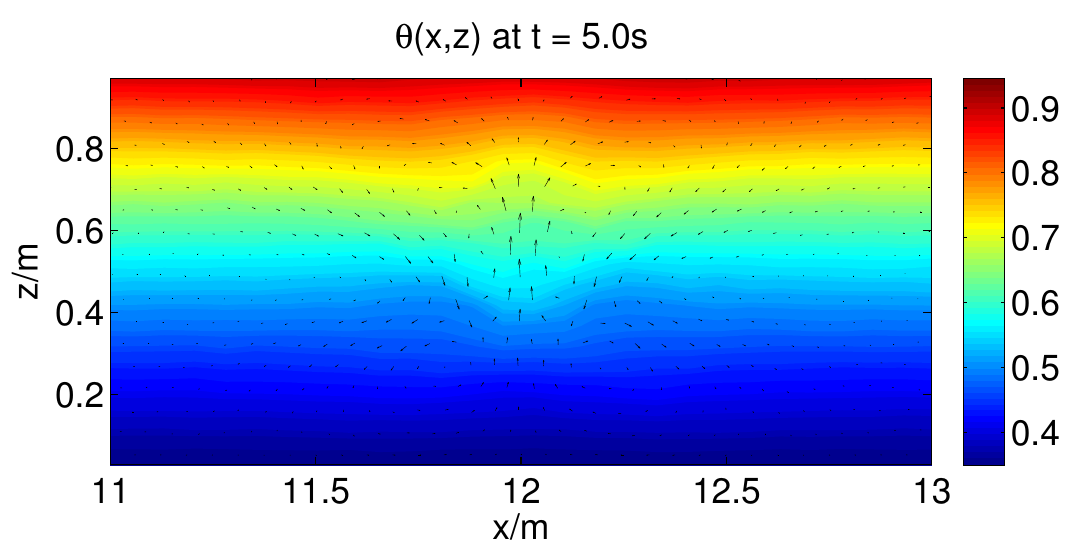}}\\ 
 \hspace*{-0.2cm}{\includegraphics[width=3in]{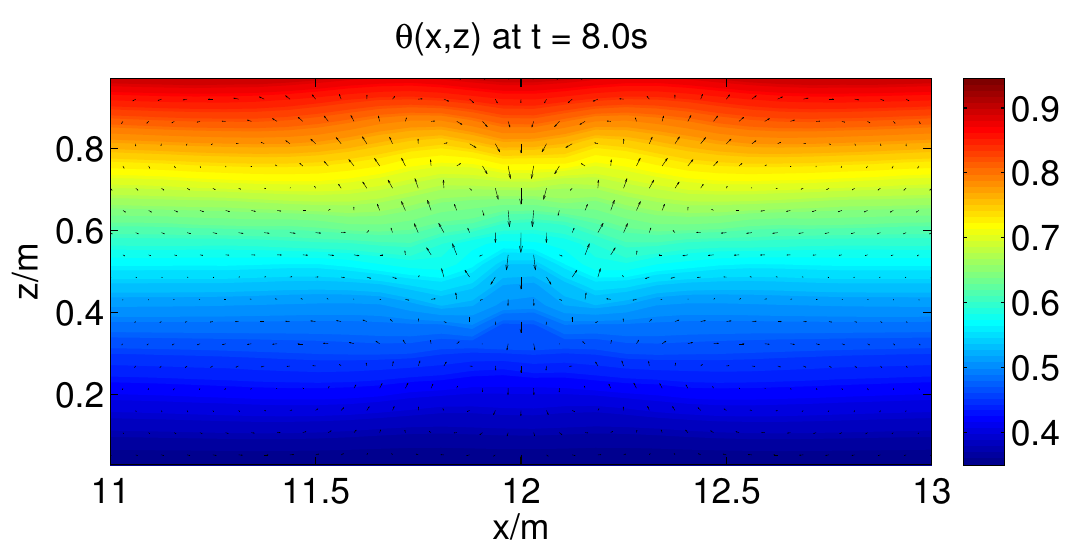}} & 
 \hspace*{-0.0cm}{\includegraphics[width=3in]{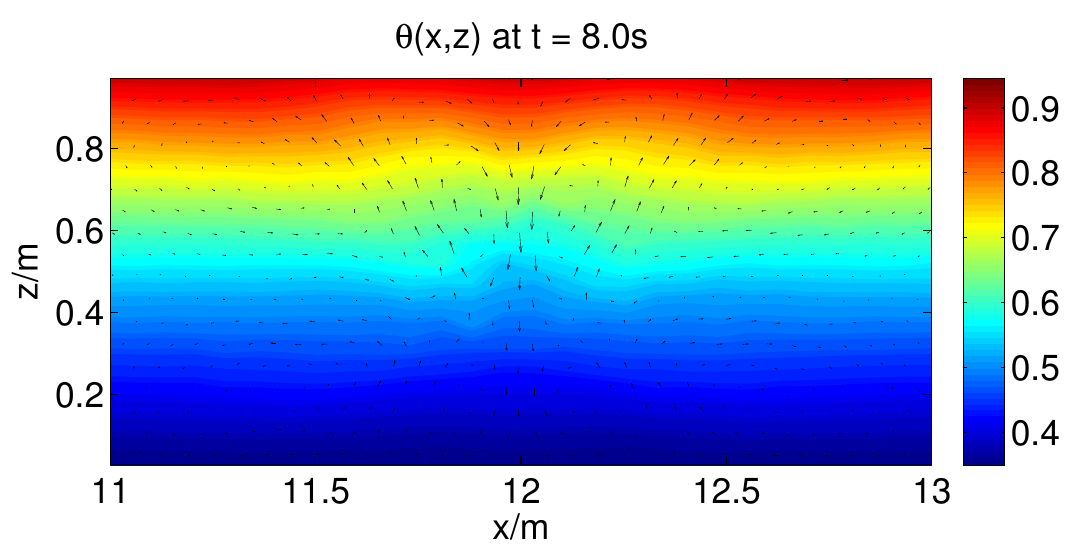}}     
 \end{tabular}
  \caption{Anelastic scheme: snapshots of the wave propagation on the regular (left column) 
  and the irregular (right column) mesh (snapshots for pseudo-incompressible scheme are very similar, 
  hence not shown)}                                                                                             
 \label{fig:anelastic_dynamics}
 \end{figure}

 Analogously, we show in Figure~\ref{fig:anelastic_dynamics} snapshots of the potential temperature 
 $\theta(x,z,t)$ of the anelastic scheme. The snapshots for the pseudo-incompressible scheme are very similar and hence not shown.
 Comparing with Figure~\ref{fig:byncy_dynamics}, the wave structure on the velocity and potential temperature fields 
 are rather similar, for both time instances and both mesh types, to those 
 obtained with the Boussinesq scheme, noticing that the magnitude of displacement of $\theta$ from equilibrium is more 
 enhanced in the anelastic and pseudo-incompressible case. Again, the irregular mesh (right column) triggers solutions that 
 are slightly non axis-symmetric with respect to $x=10\,$m, but agree in general very well with the 
 internal gravity wave propagations obtained on the regular mesh.

%  In Figure~\ref{fig:pseudo_dynamics}, we show snapshots of the potential temperature 
%  $\theta (x,z,t)$ at time $t = 5\,$s and $t = 8\,$s for a region $[11\,{\rm m},13 \,{\rm m}] \times  [0,1\,{\rm m}]$ 
%  and for the density profile (i). They show for both time instances and both mesh types very 
%  similar dynamical behavior as for Boussinesq and anelastic cases.  
%  \begin{figure}[t] \centering
%  \begin{tabular}{cc}    
%  \hspace{-0.2cm}{\includegraphics[width=3in]{23newfig_fqc_reg_pseudoincomp_2d0_bdr_t5s.pdf}} & 
%  \hspace{-0.0cm}{\includegraphics[width=3in]{25newfig_fqc_irreg_pseudoincomp_2d0_bdr_t5s.pdf}}\\ 
%  \hspace{-0.2cm}{\includegraphics[width=3in]{24newfig_fqc_reg_pseudoincomp_2d0_bdr_t8s.pdf}} & 
%  \hspace{-0.0cm}{\includegraphics[width=3in]{26newfig_fqc_irreg_pseudoincomp_2d0_bdr_t8s.pdf}}
%  \end{tabular}
%   \caption{Pseudo-incompressible scheme: snapshots of the wave propagation on regular (left column) 
%            and irregular (right column) meshes for $\bar \rho (z)= e^{-z}$.
%            }                                                                                             
%   \label{fig:pseudo_dynamics}
%  \end{figure}

\begin{figure}[t]\centering
 \begin{tabular}{cc} 
  \hspace*{-0.2cm}{\includegraphics[width=3.3in]{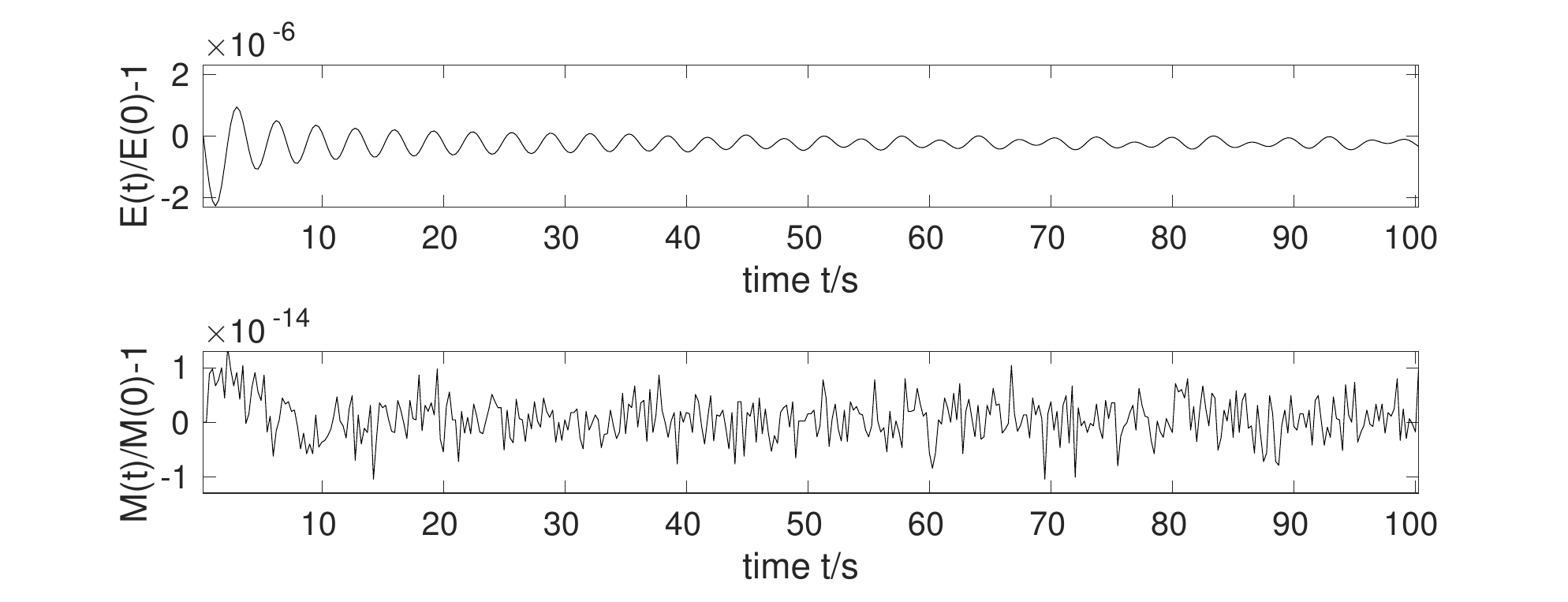}} &
  \hspace*{-0.0cm}{\includegraphics[width=3.3in]{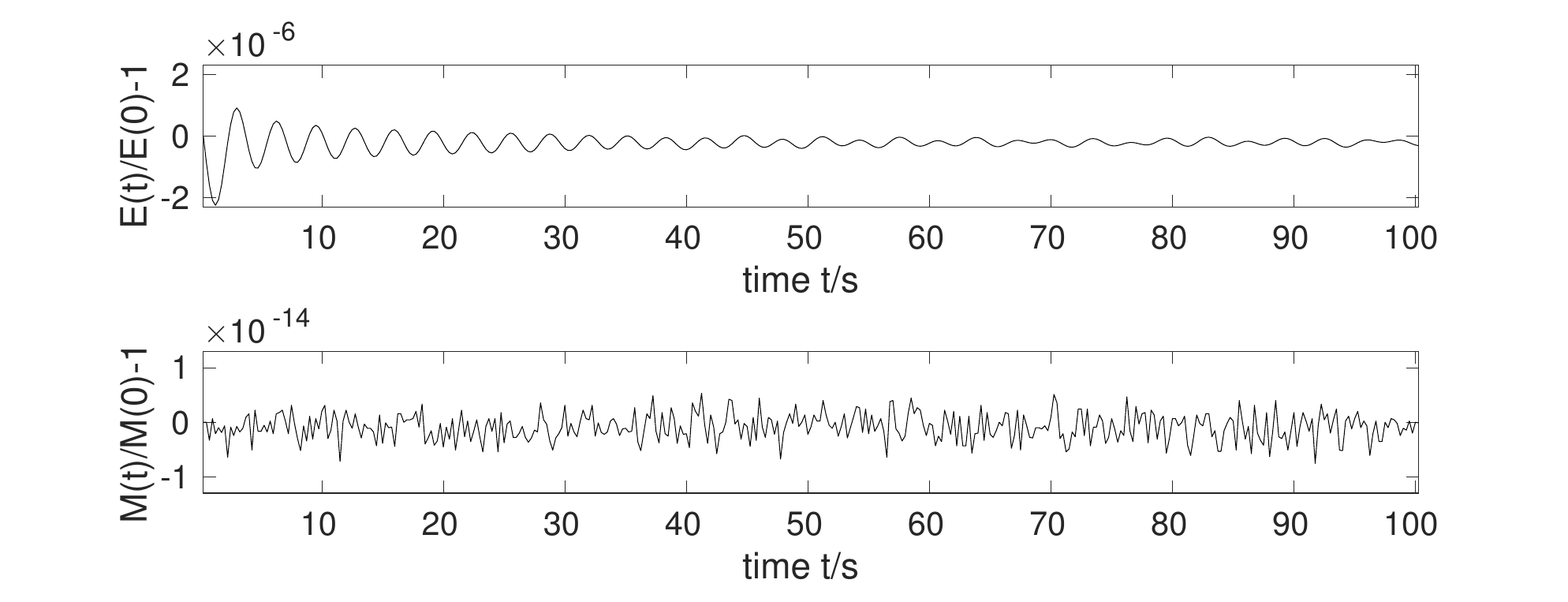}}  
  \end{tabular}
\caption{Boussinesq scheme: relative errors of total energy $E(t)$ and mass 
  $M(t)$ for the regular (left column) and the irregular (right column) mesh. 
  }                                                                                             
  \label{fig:byncy_diagnostics}
 \end{figure}

\begin{figure}[t]\centering
 \begin{tabular}{cc} 
  \hspace*{-0.2cm} {\includegraphics[width=3.3in]{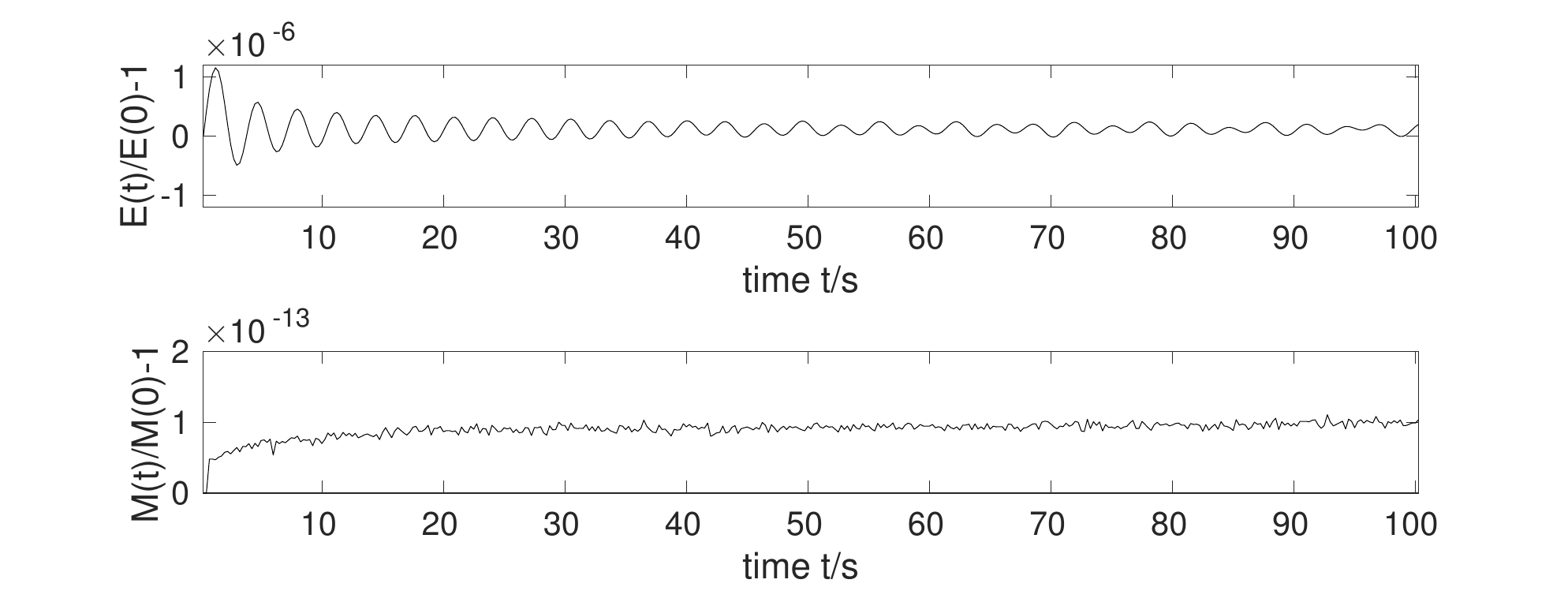}} &
  \hspace*{-0.0cm} {\includegraphics[width=3.3in]{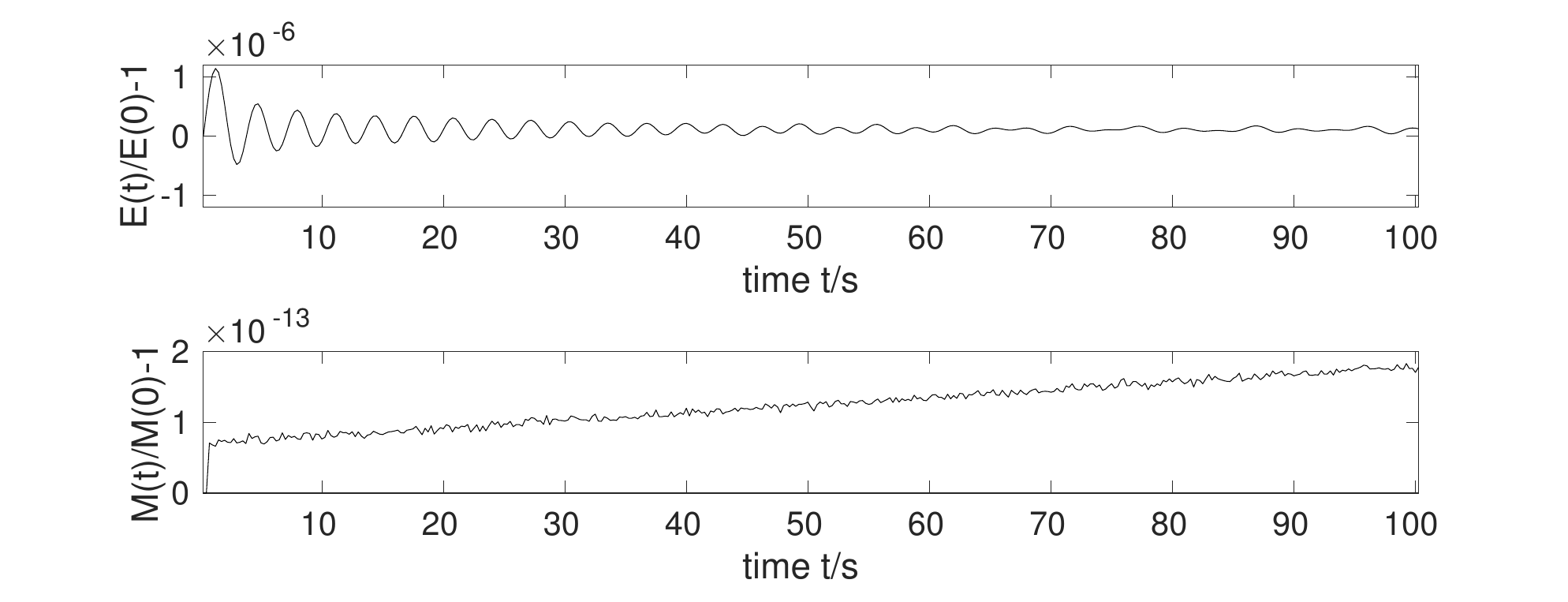}}
  \end{tabular}
\caption{Anelastic scheme: relative errors of total energy $E(t)$ and mass $M(t)$
for the regular (left column) and the irregular (right column) mesh. 
}
\label{fig:anelastic_diagnostics}
 \end{figure}

\begin{figure}[t]\centering
 \begin{tabular}{cc} 
  \hspace*{-0.2cm} {\includegraphics[width=3.3in]{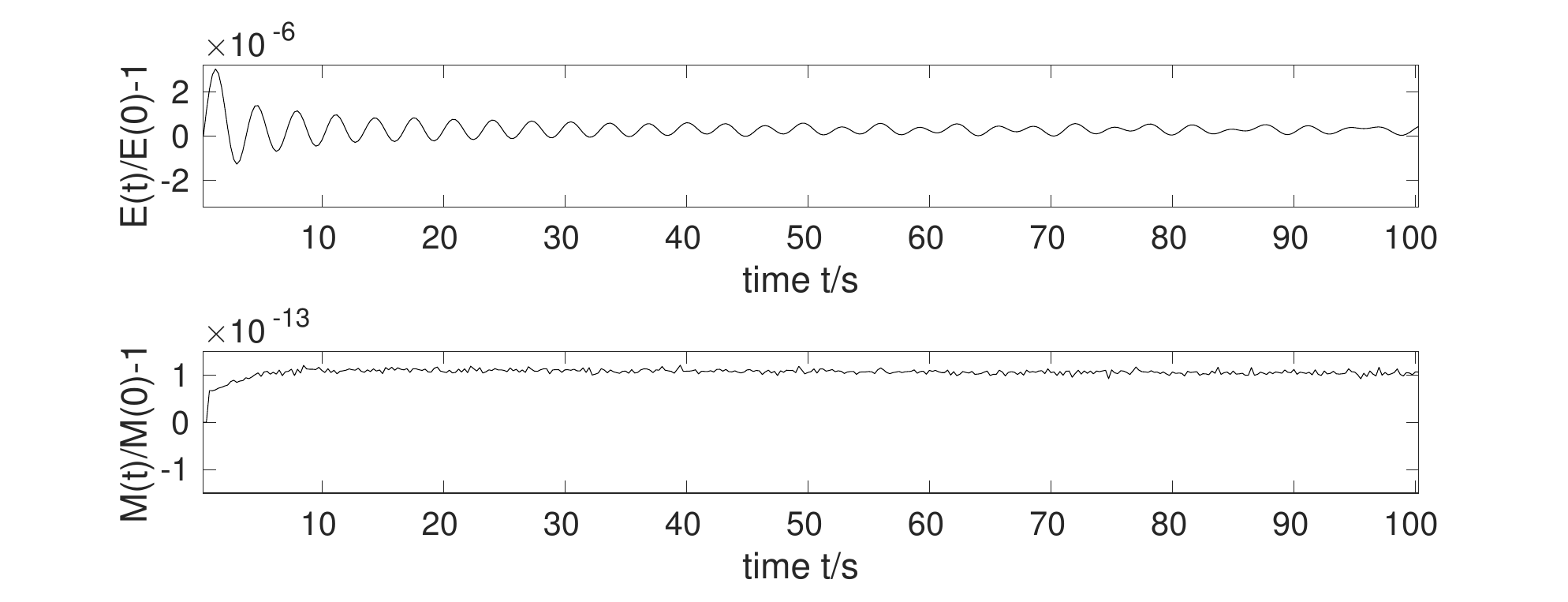}} &
  \hspace*{-0.0cm} {\includegraphics[width=3.3in]{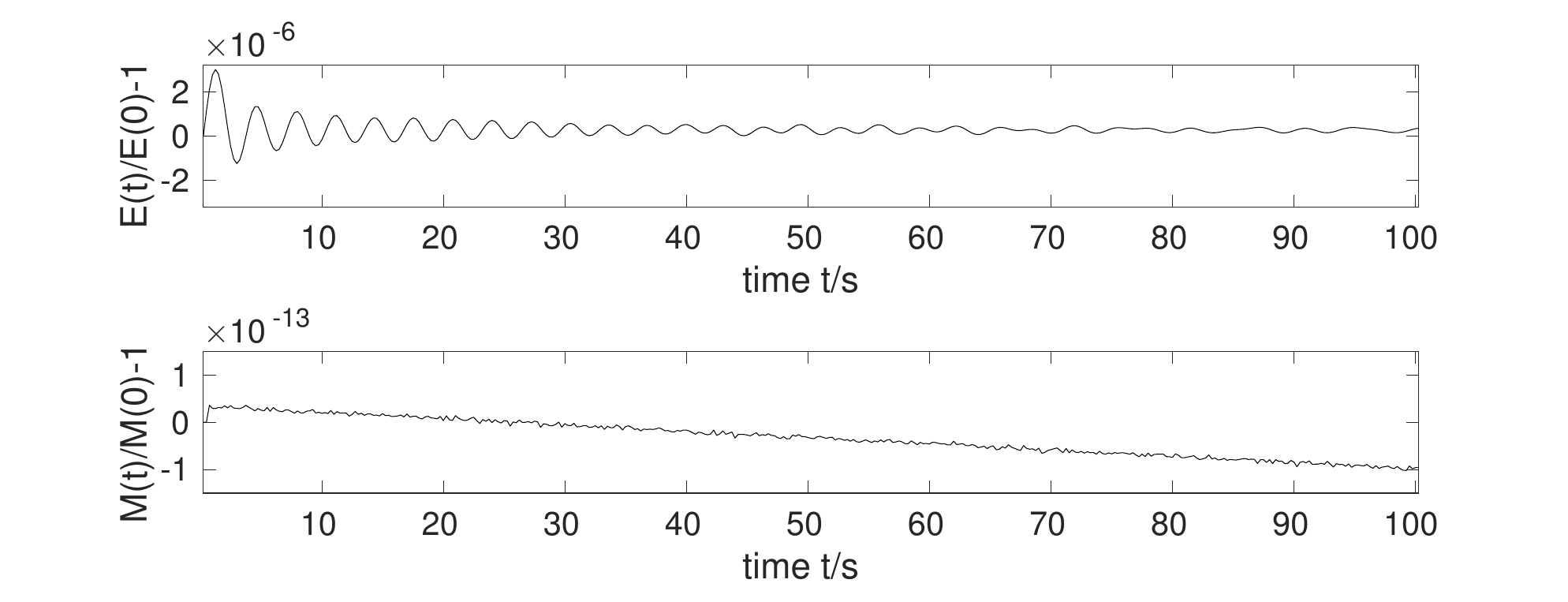}}
  \end{tabular}
\caption{Pseudo-incompressible scheme: relative errors of total energy $E(t)$ and mass 
  $M(t)$ for the regular (left column) and the irregular (right column) mesh. 
 }                                                                                             
  \label{fig:pseudo_diagnostics}
 \end{figure}

 \paragraph{Results on the conservation properties.}
 Figure~\ref{fig:byncy_diagnostics} illustrates the time evolution of the relative errors 
 (determined as ratio of current values at $t$ over initial value at $t=0$) 
 of total energy $E(t)$ (upper panels) and mass $M(t)$ %according to \eqref{equ_buyancymass_conservation} 
 (lower panels)
 of the Boussinesq scheme % equations \eqref{equ_buyancymass_conservation} 
 for the regular (left column) and the irregular (right column) mesh. 
 Analogously, Figure~\ref{fig:anelastic_diagnostics} shows these relative error 
 values for the anelastic scheme and Figure~\ref{fig:pseudo_diagnostics} for 
 the pseudo-incompressible scheme. 
   
 For all three schemes and on both mesh types, the total energy shows an oscillatory 
 behavior while being very well conserved in the mean for long integration times. The magnitudes of 
 these oscillations are at the order of $10^{-6}$, but they depend on the time step size; 
 here we used $\Delta t = 0.25\,{\rm s}$. Reducing the time step size by a factor of $10$ decreases simultaneously 
 the magnitude of the relative errors in total energy by the same factor (not shown).
 Hence, all three variational schemes show the expected $1$st-order convergence rate with time 
 (cf. time scheme derivation in \cite{PaMuToKaMaDe2010}). 
  
 In case of the Boussinesq scheme, mass is conserved at the order of $10^{-14}$
 for both the regular and the irregular mesh. 
 In the anelastic and pseudo-incompressible cases, mass is conserved at the order of $10^{-13}$ for 
 both mesh types. On the irregular mesh though we observe 
 a slight growth in the anelastic, and a slight decline in the pseudo-incompressible case,
 but within the order of $10^{-13}$ on a very acceptable level.

%  \begin{remark}\label{remark_bsnqISan}{\rm 
%  By setting $\bar \rho =\rho_0$ and $\bar \theta = \theta_0 $ as constants, 
%  the analytical Boussinesq and anelastic models agree. Similarly, by initializing 
%  the discrete schemes accordingly, i.e. by setting $\theta = -b$, $\bar \rho = 1$, $\bar\theta = \theta_0 = 1$, and $c_p = 1$,
%  we exactly recover the numerical results of the Boussinesq by the anelastic scheme 
%  (not shown).}
%  \end{remark}

 \paragraph{Investigation of the frequency representation.}
 
  We study the frequency spectra of the occurring internal gravity waves for all three schemes.
  Consider the Boussinesq system in hydrostatic equilibrium with a reference buoyancy 
  $\bar b(z)$ and a pressure $P_b$ balance like $\frac{\partial P_b}{\partial z} = \bar b$. 
  When out of equilibrium, the system tends to a hydrostatic balance 
  by emitting internal gravity waves that obey the dispersion relation
  \begin{equation}\label{equ_analytic_spectra}
   \omega^2 = \frac{k_x^2N_b^2}{\bf k^2}  
  \end{equation}
  with wave vector ${\bf k} = (k_x,k_y) \in \mathbbm Z \setminus 0$, %${\bf k} \neq 0$, 
  in which $N_b^2 :=  \frac{d \bar b}{dz}= \frac{g}{\theta_0} \frac{ d \bar \theta}{d z}$, assumed to be a constant, 
  denotes the Brunt-V\"as\"al\"a frequency
  for the case of Boussinesq equations. 
 
  For the anelastic equations, we assume
  that the reference states $\bar \rho (z)$ and $ \bar\theta (z)$ are such that
  \begin{equation}\label{N_sigma}  
   N^2 = \frac{g}{\bar \theta }\frac{ d\bar\theta }{ dz} \quad\text{and}\quad 
  \sigma_{\rm a} = \frac{1}{4} \left( \frac{1}{\bar \rho } \frac{d  \bar \rho}{dz} \right) ^2 - \frac{1}{2}\frac{d}{dz} \left( \frac{1}{\bar \rho } \frac{d \bar \rho}{dz} \right)
  \end{equation} 
  are constant numbers. Then, the dispersion relation takes the simple form
  \begin{equation}\label{equ_disprel_anelastic}
   \omega ^2 = \frac{N ^2 k _x ^2 }{ \mathbf{k} ^2  + \sigma_a  }.
  \end{equation} 
  Constant values for $N$ and $ \sigma_a $ are obtained by  
  taking $ \bar \theta (z)= \alpha e^{\frac{N ^2 }{g}z}$ and
  $ \bar \rho (z)= \beta e^{Kz}$, in which case $ \sigma_a = \frac{K ^2 }{4}$. 
  
  Similarly for the pseudo-incompressible equations in hydrostatic equilibrium,
  we assume that the reference states $\bar \rho (z)$ and $ \bar\theta (z)$ are such that
 \begin{equation}\label{N_sigma_pi}  
 N^2 = \frac{g}{\bar \theta }\frac{ d\bar\theta }{ dz} \quad\text{and}\quad 
\sigma_{\rm pi} =  \left( \frac{1}{\bar \theta  } \frac{d  \bar \theta }{dz} + \frac{1}{2\bar \rho } \frac{d  \bar \rho}{dz} \right) ^2 
- \frac{d}{dz} \left( \frac{1}{\bar \theta  } \frac{d  \bar \theta }{dz}+\frac{1}{2\bar \rho } \frac{d \bar \rho}{dz} \right) 
 \end{equation} 
 are constant numbers. The dispersion relation takes the simple form
 \begin{equation}\label{equ_disprel_pseudo}
  \omega^2 = \frac{N^2 k_x^2}{{\bf k}^2 +\sigma _{\rm pi}  } .
 \end{equation}
 
 For all three models, one observes that the frequency spectra of the internal gravity waves are 
 anisotropic and bound from above by $N_b$, respectively $N$, in the Boussinesq, respectively anelastic or 
 pseudo-incompressible case.
 To see this, consider the extremes of \eqref{equ_disprel_anelastic}, for instance, but the same reasoning 
 works for the other cases too: the lower bound at $\min(\omega) = 0$
 results from $k_x = 0$ for any $k_y > 0$ or $\sigma_a$, while the upper bound $\max(\omega) = N$  
 from $k_x \gg k_y,\sigma_a$.

 \paragraph{Results.}
 
\begin{figure}[t]
 \begin{tabular}{cc}
  \hspace*{-0.6cm}{\includegraphics[width = 3.6in]{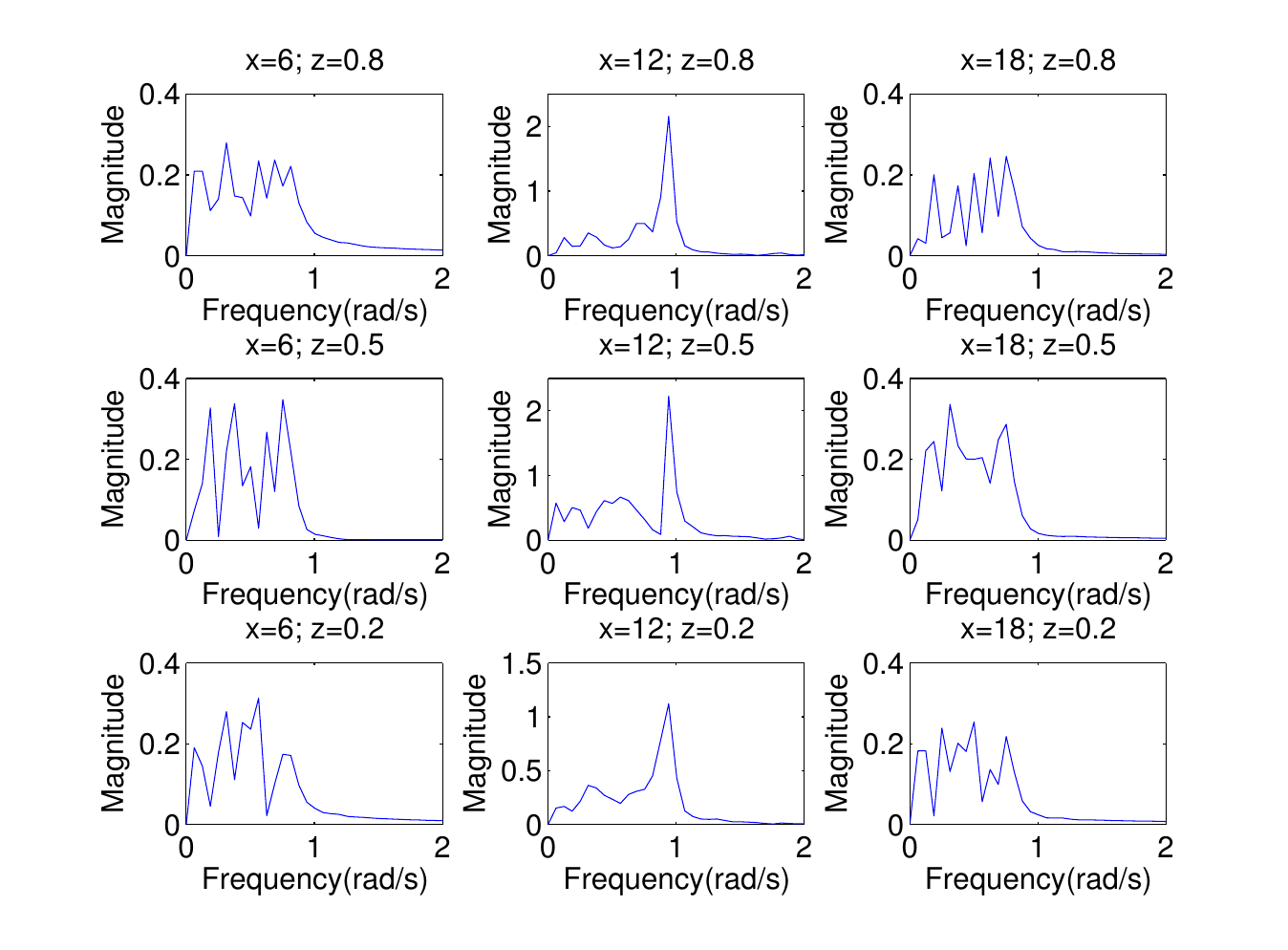}} &
  \hspace*{-0.9cm}{\includegraphics[width = 3.6in]{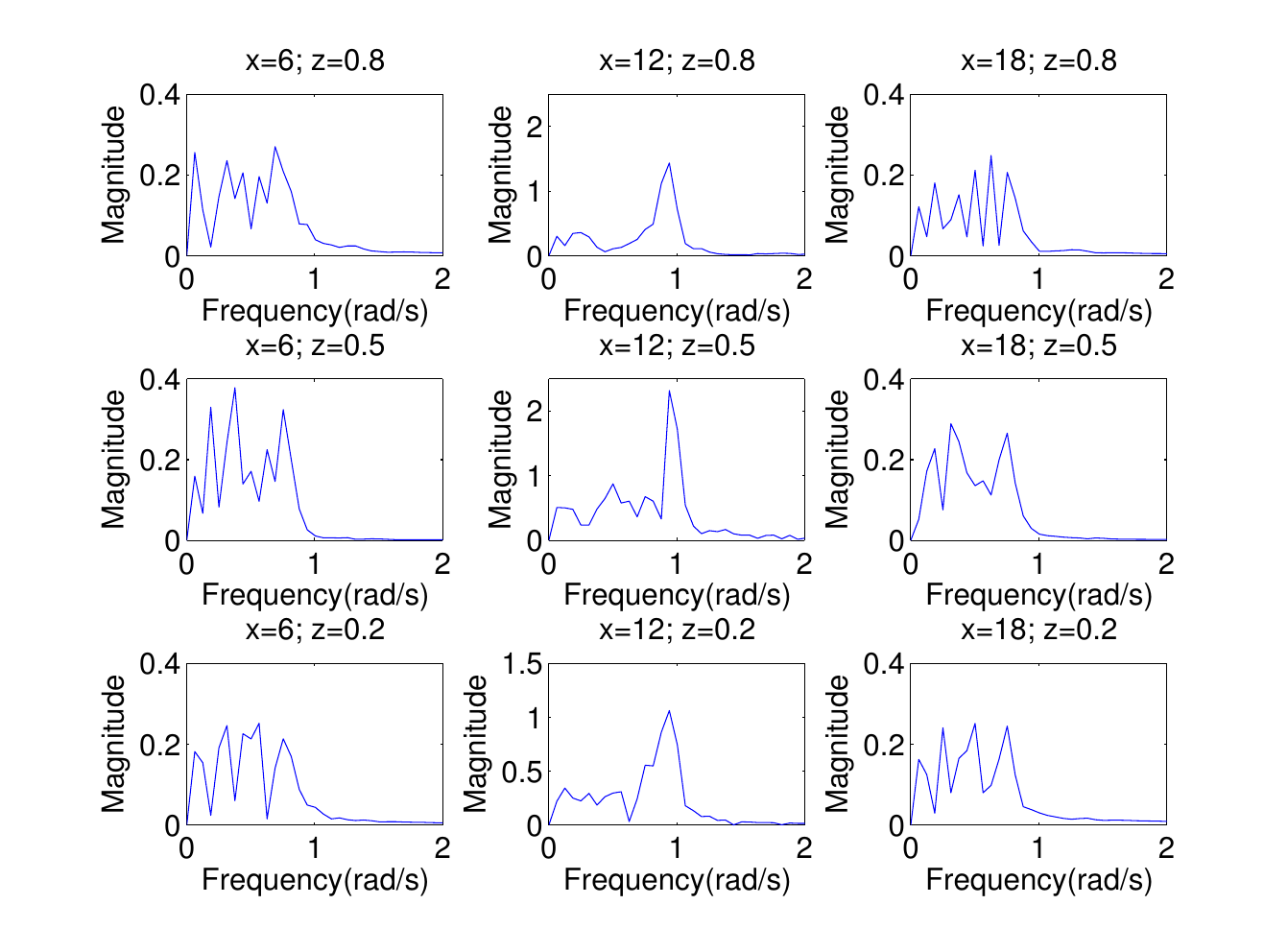}}
 \end{tabular}
  \caption{Boussinesq scheme: frequency spectra for the regular (left block) and the irregular (right block) mesh 
  determined on various points in the domain $\mathcal{D}$. The position in the panel
  indicates the corresponding position in $\mathcal{D}$, e.g. the upper left panel corresponds to a point 
  at the upper left of $\mathcal{D}$. 
  }                                                                                             
  \label{fig:byncy_freq_field}
 \end{figure}

\begin{figure}[t]
 \begin{tabular}{cc}
  \hspace*{-0.6cm} {\includegraphics[width = 3.6in]{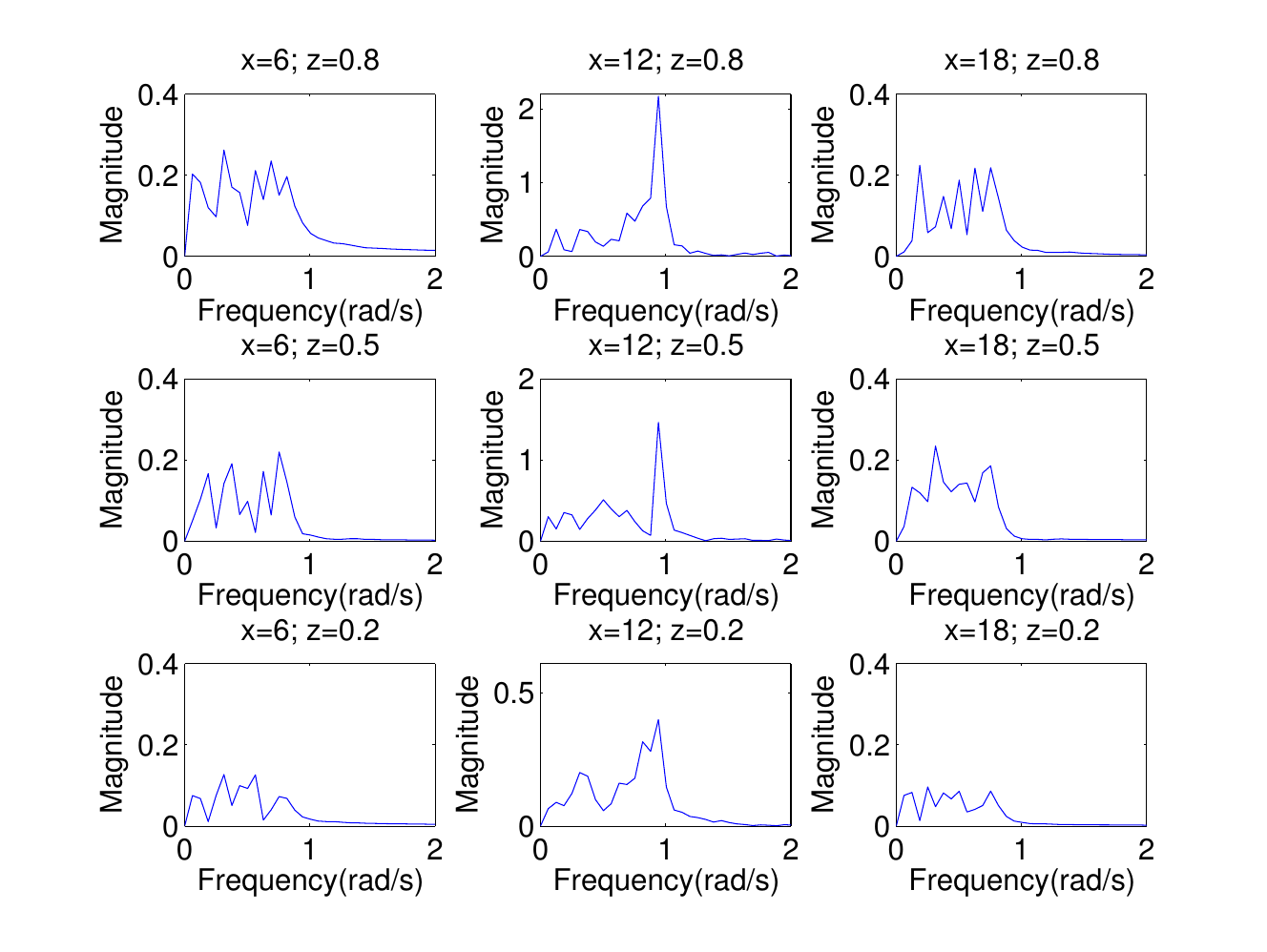}} &
  \hspace*{-0.9cm} {\includegraphics[width = 3.6in]{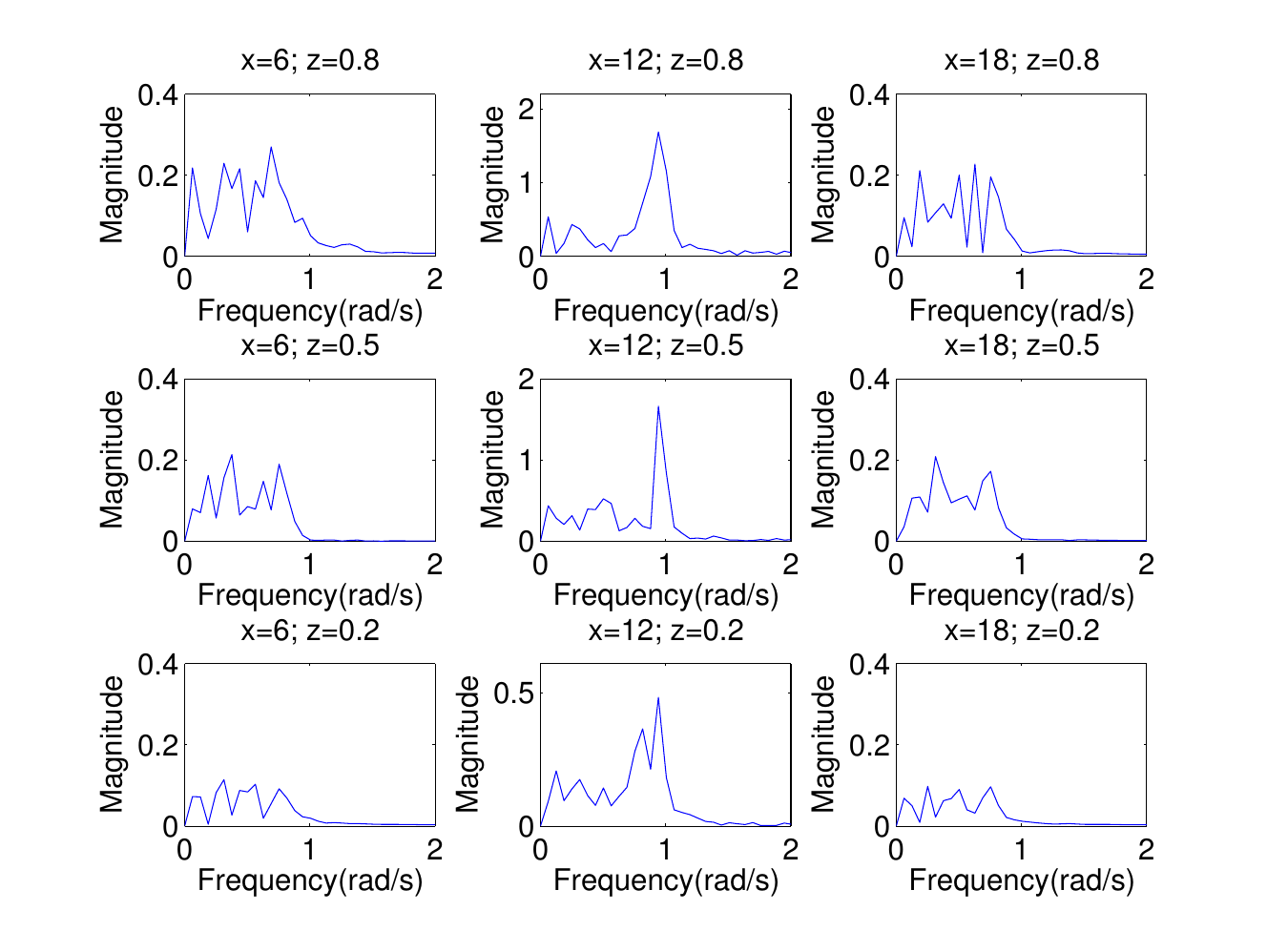}} 
 \end{tabular}  
  \caption{Anelastic scheme: frequency spectra %for initialization (i) with $\bar \rho(z) = e^{-z}$
  for the regular (left block) and the irregular (right block) mesh determined on various points in the domain $\mathcal{D}$
  similarly to Fig.~\ref{fig:byncy_freq_field}. }                                                                                             
  \label{fig:anelastic_freq_field}
 \end{figure}

\begin{figure}[t]
 \begin{tabular}{cc}
  \hspace*{-0.6cm}{\includegraphics[width=3.6in]{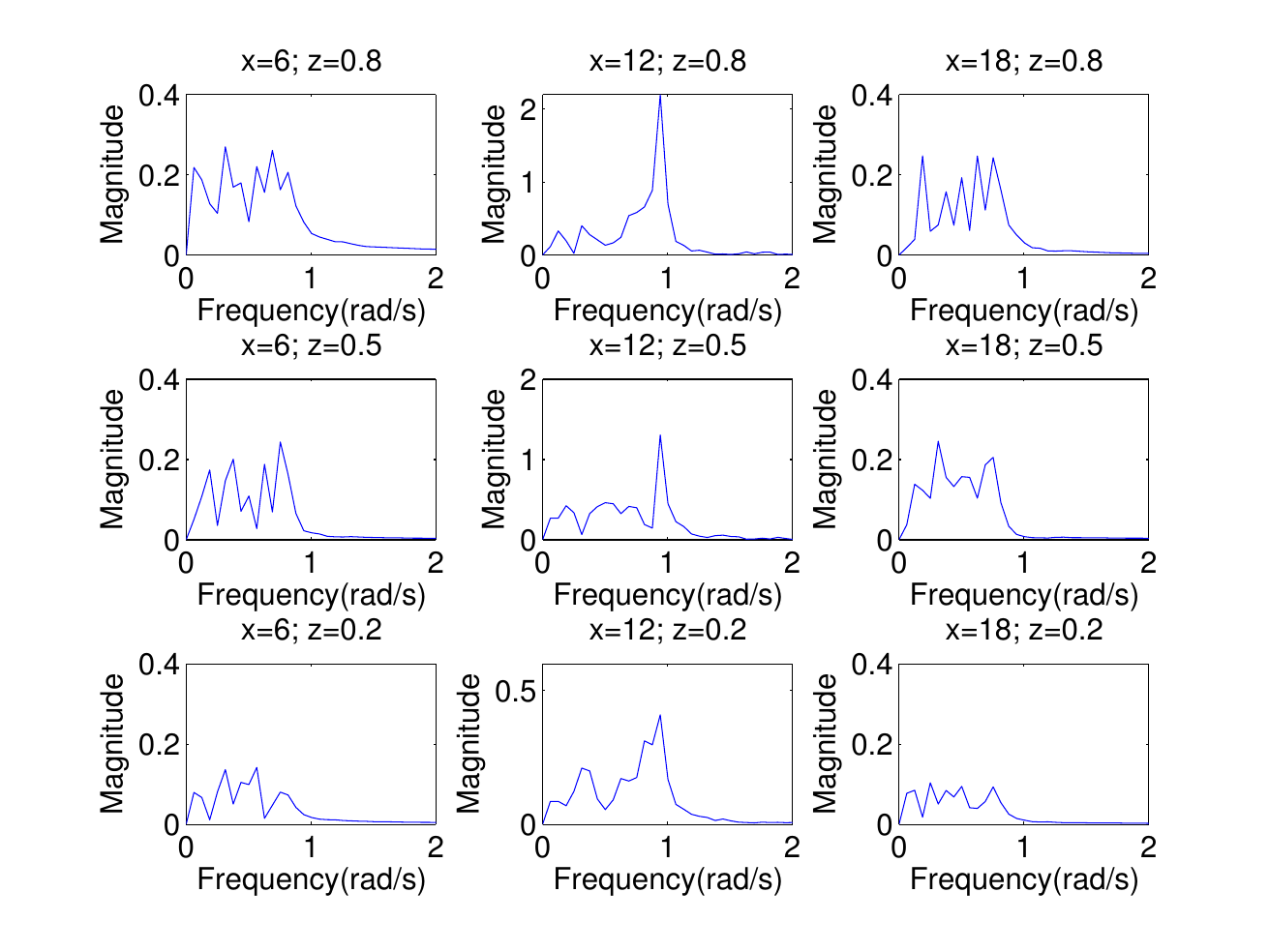}} &  
  \hspace*{-0.9cm}{\includegraphics[width=3.6in]{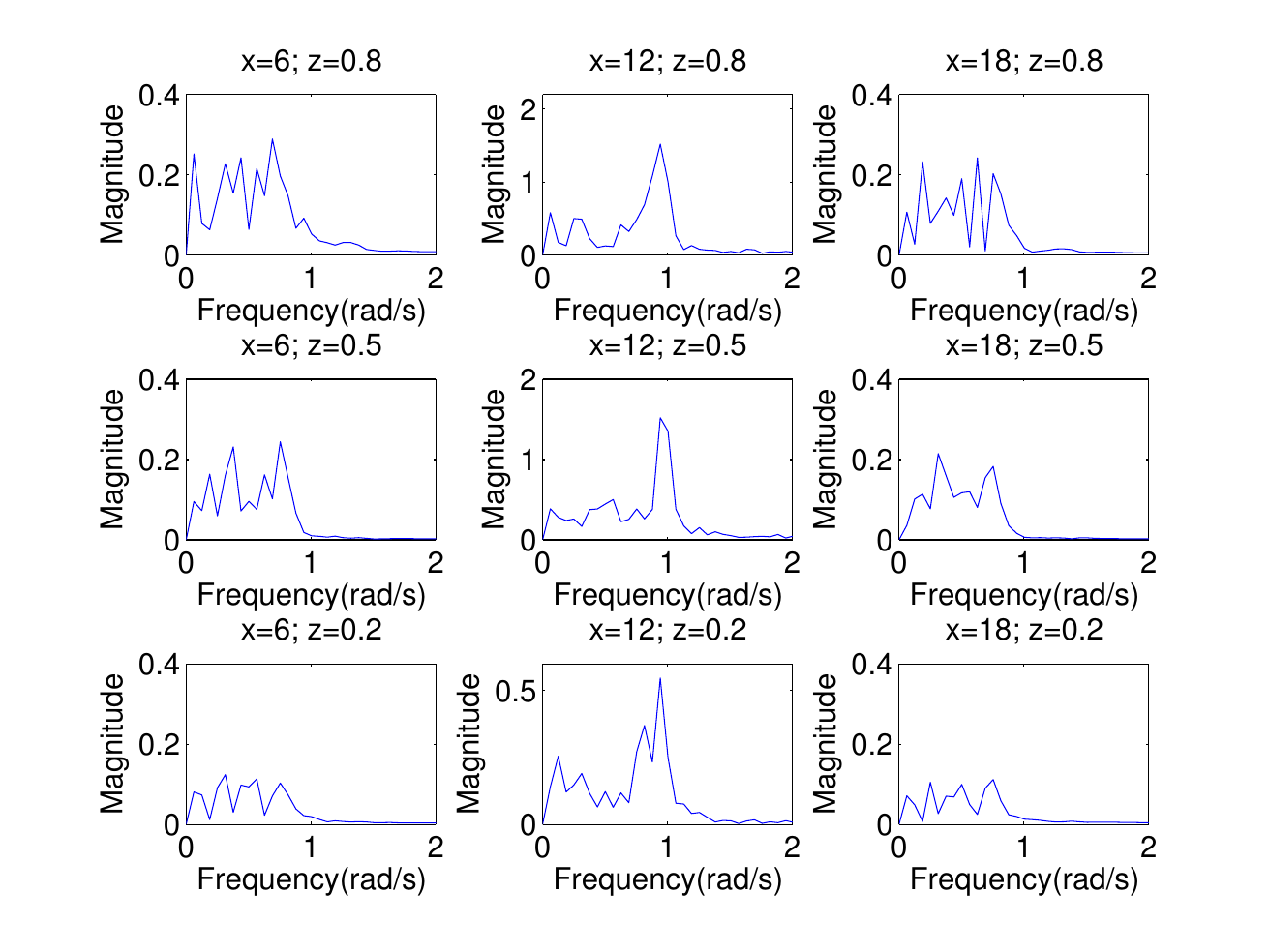}}
 \end{tabular}
  \caption{Pseudo-incompressible scheme: frequency spectra %for initialization (i) with $\bar \rho(z) = e^{-z}$
  for the regular (left block) and the irregular (right block) mesh determined on various points in the domain $\mathcal{D}  $ 
  similarly to Fig.~\ref{fig:byncy_freq_field}.  
  }                                                                                             
  \label{fig:pseudo_freq_field}
 \end{figure}

 To study numerically the dispersion relations of our discrete schemes, 
 we determine the Fourier transforms of time series of the buoyancy field $b(x,z,t)$ 
 and of the potential temperature fields $\theta(x,z,t)$ of the anelastic and pseudo-incompressible schemes
 for the time interval $t \in [0,100\,{\rm s}]$ at various locations of the computation 
 domain $\mathcal{D}$ (similar to those chosen by \cite{DeGaGBZe2014}).  
 The resulting spectra are presented in Figure~\ref{fig:byncy_freq_field}
 for the Boussinesq, Figure~\ref{fig:anelastic_freq_field} for 
 the anelastic, and Figure~\ref{fig:pseudo_freq_field} for the pseudo-incompressible schemes; 
 left blocks for the regular, and right blocks for the irregular mesh.

 For all selected sample points, these spectra show an anisotropy manifested by the 
 fact that the frequencies lie between zero and $\max (\omega) = N_b = N = 1/{\rm s}$
 with a sharp drop in the spectra at the maximal frequencies $\max (\omega)$. 
 Hence, the spectra are bound from above by $\max (\omega)$ as theoretically expected.
 Considering the central panel of each block,
 the spectrum is pronounced for values of $N_b$ in agreement with \eqref{equ_analytic_spectra}
 or of $N$ in agreement with \eqref{equ_disprel_anelastic} or \eqref{equ_disprel_pseudo}: 
 for waves with frequency near $N_b$ or $N$, the group velocity tends to zero leaving the corresponding waves trapped
 in the center of the domain. A very similar distribution of frequency spectra within the domain $\mathcal{D}$
 has been found by \cite{DeGaGBZe2014}. 
 The simulations on the irregular mesh give very similar frequency spectra.
 Hence, for all cases the spectra reflect very well the 
 properties of the analytical dispersion relations.

 \color{black}

 \section{Conclusion}
 \label{sec_conclusions}

In this paper we derived variational integrators for the anelastic and pseudo-incompressible models 
by exploiting the variational discretization framework introduced in \cite{PaMuToKaMaDe2010} for the discretization 
of incompressible fluids. In order to enable the use of this framework, we first described the anelastic and 
pseudo-incompressible approximations of the Euler equations of a perfect gas in terms of the Euler-Poincar\'e 
variational method. Applying the idea of weighted volume forms, i.e. weighted in terms of the background 
stratifications of density (anelastic) or of density times potential temperature (pseudo-incompressible), we 
could identify the appropriate groups of diffeomorphisms for the two models.

Based on these results, we defined suitable discrete versions of these diffeomorphism groups that incorporate 
the idea of weighted meshes as discrete counterparts of the weighted volume forms, in order to match the 
divergence-free conditions of the corresponding weighted velocity fields. Alongside, we defined appropriate 
weighted pairings required to derive the functional derivatives of the discrete Lagrangian that leads 
to the corresponding discrete equations of motion for the anelastic and pseudo-incompressible models, valid 
on any mesh discretization of the fluid domain. 
We then considered in detail the case of irregular 2D simplicial meshes for these two models. 
For completeness, we also considered the case of the Boussinesq equations on irregular 2D simplicial meshes, 
thereby extending the results of \cite{DeGaGBZe2014}.
For each case, we discussed the form of the equations that appears in discrete form, which is not the standard 
form in which these equations are usually written, see Table \ref{table}.

We then tested the obtained variational integrators for both regular and irregular triangular meshes by focusing 
on hydrostatic adjustment processes. These preliminary tests showed that our variational integrators capture very 
well the characteristics of the corresponding dispersion relation, in particular the upper and lower bounds of 
permitted wave numbers.  In all cases studied, both mass and energy are conserved to a high degree, following 
from the structure-preserving nature of our variational integrators.

%For acknowledgements section, please don't number the section, please begin it with \section*{Acknowledgements}
\section*{Acknowledgments} 
The authors thank D. Cugnet, M. Desbrun, F. Lott, and V. Zeitlin for helpful
discussions and valuable feedback. Both authors were partially supported by the ANR project
GEOMFLUID, ANR-14-CE23-0002-01. W. Bauer has received funding from the European Union's Horizon 2020
research and innovation programme under the Marie Sk\l odowska-Curie grant agreement No 657016.

% You may incorporate your references as follows in your main tex file.
% Using BibTex is not recommended but can be handled.

\end{document}